\documentclass{article}%

\usepackage[a4paper, margin=2.4cm]{geometry}

\usepackage{lmodern}
\usepackage[T1]{fontenc}

\usepackage{float}
\usepackage{graphicx}
\usepackage{authblk}

\usepackage[style=numeric, sorting=none]{biblatex}
\bibliography{newmark_bea}
\usepackage{multirow}%
\usepackage{amsmath,amssymb,amsfonts}%
\usepackage{amsthm}%
\usepackage{mathrsfs}%
\usepackage[title]{appendix}%
\usepackage{xcolor}%
\usepackage{bm}
\usepackage{physics}
\usepackage{xspace}

\usepackage{listings}
\usepackage{hyperref}

\newcommand{\leftf}{\mathopen{}\mathclose\bgroup\left}
\newcommand{\rightf}{\aftergroup\egroup\right}

\newcommand{\bigleftf}{\mathopen{}\mathclose\bgroup\bigl}
\newcommand{\bigrightf}{\aftergroup\egroup\bigr}

\newcommand{\ctilde}{\tilde}
\newcommand{\chat}{\widehat}

\newcommand{\ve}[1]{\mathbf{\bm{#1}}}
\newcommand{\mx}[1]{\mathbf{\bm{#1}}}
\newcommand{\yv}{\ve{y}}
\newcommand{\yvj}{\ve{y}^j}
\newcommand{\yvnull}{\ve{y_0}}
\newcommand{\yvjp}{\ve{y}^{j+1}}
\newcommand{\tj}{t^j}

\newcommand{\yvdot}{\dot{\ve{y}}}
\newcommand{\fv}{\ve{f}}
\newcommand{\RR}{\mathbb{R}}
\newcommand{\Dt}{\Delta t}

\newcommand{\fvt}{\ve{\ctilde{f}}}
\newcommand{\yvt}{\ve{\ctilde{y}}}
\newcommand{\yvtdot}{\dot{\ve{\ctilde{y}}}}

\newcommand{\xt}{\ctilde{x}}
\newcommand{\vt}{\ctilde{v}}
\newcommand{\vtdot}{\dot{\ctilde{v}}}
\newcommand{\xtdot}{\dot{\ctilde{x}}}
\newcommand{\xtddot}{\ddot{\ctilde{x}}}

\newcommand{\D}{\mx{D}}
\newcommand{\atyvtt}{\bigleftf( \yvt \leftf( t \rightf) \bigrightf)}

\newcommand{\Ord}{\mathcal{O}}

\newcommand{\Fv}{\ve{F}}
\newcommand{\qv}{\ve{q}}
\newcommand{\vv}{\ve{v}}
\newcommand{\av}{\ve{a}}

\newcommand{\qvj}{\ve{q}^j}
\newcommand{\vvj}{\ve{v}^j}
\newcommand{\avj}{\ve{a}^j}
\newcommand{\Fvjp}{\ve{F}^{j+1}}
\newcommand{\qvjp}{\ve{q}^{j+1}}
\newcommand{\vvjp}{\ve{v}^{j+1}}
\newcommand{\avjp}{\ve{a}^{j+1}}

\newcommand{\Mm}{\mx{M}}
\newcommand{\Cm}{\mx{C}}
\newcommand{\Km}{\mx{K}}
\newcommand{\Mminv}{\mx{M}^{-1}}

\newcommand{\Ctm}{\mx{\ctilde{C}}}
\newcommand{\Ktm}{\mx{\ctilde{K}}}
\newcommand{\Fvft}{\Fv \leftf( t \rightf)}
\newcommand{\Fvftp}{\Fv' \leftf( t \rightf)}
\newcommand{\Fvftpp}{\Fv'' \leftf( t \rightf)}
\newcommand{\Ftvft}{\ve{\ctilde{F}} \leftf( t \rightf)}

\newcommand{\Nm}{\mx{0}}
\newcommand{\Nv}{\ve{0}}
\newcommand{\Em}{\mx{I}}

\newcommand{\Gm}{\mx{G}}
\newcommand{\Hm}{\mx{H}}

\newcommand{\Amtqv}{\mx{A}\leftf(\tau, \qv,\vv\rightf)}
\newcommand{\Am}{\mx{A}}

\newcommand{\Bmfhg}{\mx{B}\leftf(\Dt, \gamma \rightf)}

\newcommand{\Fvftau}{\Fv \leftf( \tau \rightf)}
\newcommand{\Fvftaup}{\Fv' \leftf( \tau \rightf)}
\newcommand{\Fvftaupp}{\Fv'' \leftf( \tau \rightf)}

\newcommand{\qvnull}{\ve{q}_0}
\newcommand{\qdotvnull}{\ve{v}_0}

\newcommand{\Chm}{\mx{\chat{C}}}
\newcommand{\Khm}{\mx{\chat{K}}}
\newcommand{\Fhvft}{\chat{\Fv} \leftf( t \rightf)}
\newcommand{\Fhvftp}{\chat{\Fv}' \leftf( t \rightf)}
\newcommand{\Fhvftpp}{\chat{\Fv}'' \leftf( t \rightf)}

\newcommand{\Fhvtwoft}{\chat{\Fv}_2 \leftf( t \rightf)}

\newcommand{\tstar}{t^{*}}

\newcommand{\sign}{\operatorname{sign}}

\newcommand{\compensation}{compensation\xspace }
\newcommand{\compensations}{compensations\xspace }
\newcommand{\compensated}{compensated\xspace }

\newcommand{\compensating}{compensating\xspace }



\begin{document}

\title{Improving the accuracy of the Newmark method through backward error analysis}

\author[1,2,*]{Don\'{a}t M. Tak\'{a}cs}

\author[1,2]{Tam\'{a}s F\"{u}l\"{o}p}

\affil[1]{Department of Energy Engineering, Faculty of Mechanical Engineering, Budapest University of Technology and Economics, M\H{u}egyetem rkp.~3., Budapest, H-1111, Hungary}

\affil[2]{Montavid Thermodynamic Research Group, c/o ETTE, Lovas \'{u}t~18., Budapest, H-1012, Hungary}

\affil[*]{Corresponding author, e-mail: takacs@energia.bme.hu}

\maketitle

\begin{abstract}
We use backward error analysis for differential equations to obtain modified or distorted equations describing the behaviour of the Newmark scheme applied to the transient structural dynamics equation. Based on the newly derived distorted equations, we give expressions for the numerically or algorithmically distorted stiffness and damping matrices of a system simulated using the Newmark scheme. Using these results, we show how to construct compensation terms from the original parameters of the system, which improve the performance of Newmark simulations. The required compensation terms turn out to be slight modifications to the original system parameters (e.g. the damping or stiffness matrices), and can be applied without changing the time step or modifying the scheme itself. Two such compensations are given: one eliminates numerical damping, while the other achieves fourth-order accurate calculations using the traditionally second-order Newmark method. The performance of both compensation methods is evaluated numerically to demonstrate their validity, and they are compared to the uncompensated Newmark method, the generalized-$\alpha$ method and the 4th-order Runge--Kutta scheme.
\end{abstract}

\section{Introduction}\label{sec:introduction}

The Newmark scheme has been extensively used for numerically solving structural mechanics problems since its inception \cite{newmark1959method} in the late 1950s. Along with its various extensions and generalizations (e.g.~\cite{wood1980alpha,chung1993time,hilber1997improved}), it is still widely used for predicting the temporal behaviour of various systems under external loads \cite{bathe2014finite,hughes2012finite,bathe2014frontiers}. Even though the state of the art in time integration of finite element models has advanced beyond these methods \cite{brun2015two,malakiyeh2018development,noh2019direct,bamer2021newmark,sanchez2021symplectic,stasio2021explicit,zakian2021transient,choi2022time,huang2022spacetime,cockburn2023combining,dvorak2023energy,huang2023HHT,soares2023enhanced,luo2024nonsmooth}, most commercially available and widely used finite element software packages still rely primarily on the Newmark method or its extension, the generalized HHT-$\alpha$ method for solving mechanical problems \cite{ansys2023implicit,abaqus2023implicit}.

For describing the behaviour of geometric or structure-preserving numerical schemes applied to systems of ordinary differential equations (ODEs), the mathematical technique of backward error analysis (BEA) has been developed (see e.g.~\cite{griffiths1986scope,reich1999backward,hairer2000asymptotic}, among others). Application of this approach to a set of ODEs and a numerical method yields the so-called modified or distorted equation, which describes the behaviour of the discrete-time numerical method as a system of continuous-time ODEs. Hence, distorted equations produced by BEA allow for the qualitative and quantitative analysis of numerical methods, as they can be compared to the original equations directly, using the same mathematical tools.

While there have been results connecting the Newmark method and structure-preserving schemes \cite{simo1992exact,kane2000variational}, the application of backward error analysis to obtain distorted equations corresponding to the Newmark scheme -- to the best of our knowledge -- has not yet been published. The present paper bridges this gap: we give two forms of the distorted equation for linear systems under transient excitations simulated using the Newmark scheme. One is of a first-order form, treating position and velocity as separate variables; the other is a second-order form which can directly be compared to the original system. (Two partially related, but distinct approaches are introduced in \cite{wood1986unified,krenk2006energy}.)

% \cite{wood1986unified} backward error analysist emleget, de az valami más dolog, és csak a merevséget és a csillapítást módosítja.

Using the results outlined above, we then introduce two constructions for compensating qualitative and numerical errors introduced by the Newmark method. First, we show how the numerical damping introduced by the Newmark method can be mitigated for arbitrary values of the Newmark parameters $\gamma$ and $\beta$. Second, we show how the Newmark scheme (which is traditionally second-order accurate at most) can be used to perform fourth-order accurate calculations for $\gamma=1/2$ and $\beta=1/6$. 

Both constructions use backward error analysis-based compensation, a technique introduced in this paper. During compensation, the original numerical method is not modified at all. However, by knowing the distortions introduced by the application of the numerical scheme from BEA, the parameters of the original system (such as its damping, stiffness or excitation) can be changed slightly so that the undesirable effects of the numerical method are cancelled out. (This approach is somewhat reminiscent of, but unrelated to, input shaping \cite{robinett2002input} in robotics.)

The main advantage of the compensation technique is that the performance of the numerical method can be improved solely through the appropriate tuning of the original system parameters; thus, the improvements can be readily used in existing software. We see this as a significant advantage over introducing novel numerical methods, which often need third-party implementations to reach more widespread adoption.

The outline of the paper is as follows. First, we give an overview of the backward error analysis of ODEs in general in Section~\ref{sec:bea_illustration}. This is followed by the application of BEA to the Newmark method to obtain the distorted ODEs in Section~\ref{sec:newmark_modified}. Building on these results, we introduce the compensation technique for eliminating numerical damping and achieving fourth-order accuracy in the Newmark method in Section~\ref{sec:newmark_improvement}. Numerical examples verifying and demonstrating the results are given throughout, accompanying the respective calculations.

\section{Backward error analysis for numerical methods}\label{sec:bea_illustration}
Given a system of $n$ first-order ordinary differential equations (ODE) and an initial condition (IC):
\begin{gather}
    \yvdot\leftf(t\rightf) = \fv\leftf( \yv\leftf(t\rightf), t \rightf),\quad \yv(0)=\yvnull,  \label{eq:ode}
\end{gather}
with a solution $\yv: \RR \rightarrow \RR^n$, initial condition $\yvnull$ and non-autonomous vector field $\fv$$: \RR^n \times \RR \rightarrow \RR$, we consider a numerical method $\Phi_{\Dt}$ (with fixed time step $\Dt$) which generates numerical solutions $\yvj$ at discrete time instants $\tj := j \Dt$ for $j=0, 1, 2, \ldots J$.

It is well known that, generally, a numerical method solving \eqref{eq:ode} does not give the exact solution at the discrete time instants considered, i.e.\ $\yvj \neq \yv\leftf( \tj \rightf)$, though a \emph{consistent} method gives an accurate solution in the zero time-step limit, i.e. it fulfils the condition $\lim_{\Dt \to 0} \yvj = \yv\leftf( \tj \rightf)$. Traditional or forward error analysis considers the accuracy of an exact forward step $\yv\leftf( t+\Dt \rightf)$ versus a numerical step induced by $\Phi_{\Dt}$. On the other hand, backward error analysis (BEA) of numerical methods (as given by \cite{griffiths1986scope,reich1999backward,hairer2000asymptotic,moan2006modified}, among others) considers the existence and behaviour of the so-called \emph{modified} or \emph{distorted equation}%
\footnote{
%We will primarily use the latter term as we feel that it conveys the essence of the analysis better, and avoids confusion with the idea of \emph{\compensation} introduced later on. It is worth noting, however, that the former is a more widespread terminology in the mathematics community.
Since the
frequently used
word ''modified'' 
-- primarily originating from the mathematics community --
could refer here to both the modification \emph{caused by the numerical scheme} and the modification \emph{performed by us} to counterbalance it, throughout the paper we use ''distorted'' for the former and ''compensated'' for the latter.\label{footnote:distorted}
}
of a numerical method, which fulfils the condition
\begin{gather}
    \yvt\leftf( \tj \rightf) = \yvj, \qquad \forall j=0,1,\ldots,J \label{eq:condition}
\end{gather}
and behaves according to the ODE and initial value condition
\begin{gather}
    \yvtdot\leftf(t\rightf) = \fvt\leftf( \yvt\leftf(t\rightf), t \rightf),\quad \yvt\leftf( 0 \rightf) = \yvnull, \label{eq:modeq}
\end{gather}
with $\yvt(t)$ being the solution of the distorted equation \eqref{eq:modeq}, containing the
% modified or
distorted vector field (DVF%
\footnote{%
In line with Footnote~\ref{footnote:distorted}, we will use the initialism ''DVF'' instead of the frequently used ''MVF'', which abbreviates ''modified vector field''.%
}%
)
denoted as $\fvt$. In other words, the numerical method $\Phi_{\Dt}$ is an \emph{exact integrator} of $\eqref{eq:modeq}$.

% Ezek nem inkább introduction?

\subsection{Asymptotic expansion of the distorted vector field}
In what follows, we restrict ourselves to an autonomous vector field $\fv(\yv) \equiv \fv(\yv, t)$. (Non-autonomous systems can be brought to this form by an extension of the state space, as will be shown in Section~\ref{sec:newmark_modified}.) Thus, the ODE to be solved is
\begin{gather}
    \yvdot\leftf(t\rightf) = \fv\leftf( \yv\leftf(t\rightf) \rightf),\quad \yv\leftf( 0 \rightf) = \yvnull. \label{eq:ode_aut}
\end{gather}

There are several distinct but eventually equivalent approaches for obtaining the distorted equation as an asymptotic series of the time step $\Dt$: some notable examples are \cite{reich1999backward, gonzalez1999qualitative, hairer2000asymptotic}. Here, we follow the approach of \cite{hairer2000asymptotic,moan2006modified} for a general exposition on the construction of the distorted equation, with a slightly different logic. It will become clear later that this original approach needs to be extended for the Newmark method. For the convenience of the Reader, we also provide a higher order of expansion during the calculations regarding the power series. 

Assuming that the continuous distorted equation for this method exists with solution $\yvt(t)$, a Taylor-series expansion for a time step $\Dt$ can be given as
\begin{gather}
    \yvt\leftf( t +\Dt \rightf) = \yvt( t ) + \Dt \dv{y}{t} \leftf( t \rightf) + \frac{\Dt^2}{2!}  \dv[2]{y}{t} \leftf( t \rightf) + \frac{\Dt^3}{3!} \dv[3]{y}{t} \leftf( t \rightf) + \frac{\Dt^4}{4!} \dv[4]{y}{t} \leftf( t \rightf) + \ldots \label{eq:taylor}
\end{gather}
which can also be expressed using $\fvt(\yvt)$ by using an autonomous version of~\eqref{eq:modeq} and the chain rule. This gives
\begin{align}
    \yvt\leftf( t +\Dt \rightf) = & \: \yvt( t ) + \Dt\;\! \fvt \atyvtt  + \frac{\Dt^2}{2!} \;\! \leftf( \D\fvt \fvt \rightf) \atyvtt \mathrel + \nonumber \\ \nonumber
                                  & \mathrel + \frac{\Dt^3}{3!} \;\! \leftf( \D^2 \fvt\leftf(\fvt, \fvt\rightf) + \D\fvt \D\fvt \fvt \rightf) \atyvtt \mathrel + \\ 
    & \mathrel + \frac{\Dt^4}{4!} \;\! \leftf[
        \D^3 \fvt\leftf(\fvt, \fvt, \fvt \rightf) + 
        3 \D^2 \fvt\leftf(\D\fvt \fvt, \fvt\rightf) \mathrel +
        \rightf.
        \nonumber
        \\
    & \qquad \qquad \leftf. \mathrel + 
        \D\fvt \D^2 \leftf( \fvt, \fvt \rightf) +
        \D\fvt \D\fvt \D\fvt \fvt
    \rightf]\atyvtt
    + \ldots \label{eq:taylor2}
\end{align}
where $\D\fvt$ is the Jacobian of $\fvt$, and $\D^k\fvt$ are its $k$th-order vectorial derivatives.

The DVF is to be expressed as an asymptotic series in the form of
\begin{gather}
    \fvt(\yvt) = \fv(\yvt) + \Dt\;\! \fv_1(\yvt) + \Dt^2 \;\! \fv_2(\yvt) + \Dt^3 \;\! \fv_3(\yvt) + \ldots \label{eq:asymptotic}
\end{gather}
One fundamental reason for this is that the DVF of an autonomous vector field might be non-autonomous due to a small, periodic perturbation \cite{moan2006modified,oneale2009preservation}, thus the above series might not actually be convergent, and in such cases the left side is only a formal expression. However, truncated at a chosen power of the time step, \eqref{eq:asymptotic} is a valid approximation of the DVF up to and including that order \cite{reich1999backward, hairer2000asymptotic, hairer2006geometric}.

Substituting \eqref{eq:asymptotic} into \eqref{eq:taylor2} and collecting powers of~$\Dt$ yields
\begin{align}
    & \yvt\leftf( t +\Dt \rightf)
    \nonumber
    \\*
    & = \yvt(t) + 
    \Dt\;\! \fv \atyvtt + 
    \Dt^2\;\! \left(
    \fv_1 
    + \tfrac{1}{2} \D\fv \fv
    \right) \atyvtt 
    \mathrel +
    \nonumber
    \\* \nonumber
    & \quad + \Dt^3\;\! 
    \left\lbrace
        \fv_2 
        +
        \tfrac{1}{2}
        \left[
              \D \fv \fv_1 
            + \D \fv_1 \fv 
        \right]
        +
        \tfrac{1}{6}
        \left[
           \D^2 \fv \leftf( \fv, \fv \rightf) 
        +  \D\fv \D\fv \fv 
        \right]
    \right\rbrace 
    \atyvtt 
    \mathrel +
    \\* \nonumber
    &  \quad + 
    \Dt^4 \left\lbrace
        % Wolfram-sorrend a multilinear_Dn_expansion.nb alapján
        % \fv_3
        % + \tfrac{1}{2} \D\fv \fv_2
        % + \tfrac{1}{6} \D\fv \D\fv \fv_1
        % + \tfrac{1}{24} \D\fv \D\fv \D\fv \fv
        % + \tfrac{1}{6} \D\fv \D\fv_1\fv
        % +
        % \right. \\ \nonumber \left.
        % + \tfrac{1}{24} \D\fv \D^2\fv \leftf( \fv, \fv \rightf)
        % + \tfrac{1}{2} \D\fv_1 \fv_1
        % + \tfrac{1}{6} \D\fv_1 \D\fv \fv
        % + \tfrac{1}{2} \D\fv_2 \fv
        % + \tfrac{1}{3} \D^2\fv \leftf(\fv, \fv_1 \rightf) 
        % +
        % \right. \\ \nonumber \left.
        % + \tfrac{1}{8} \D^2\fv \leftf(\fv, \D\fv\fv \rightf)
        % + \tfrac{1}{6} \D^2\fv_1 \leftf(\fv, \fv \rightf)
        % + \tfrac{1}{24} \D^3\fv \leftf(\fv, \fv, \fv \rightf)
        \fv_3
        + \tfrac{1}{2} 
        \left[
            \D\fv_1 \fv_1
            + \D\fv \fv_2
            + \D\fv_2 \fv
        \right]
        \mathrel +
        \right. \\ 
        & \quad
        \nonumber \left.
         \mathrel + \tfrac{1}{6} 
        \left[
        2\:\! \D^2\fv \leftf(\fv, \fv_1 \rightf) 
        + \D\fv \D\fv \fv_1
        + \D\fv \D\fv_1\fv
        + \D\fv_1 \D\fv \fv
        + \D^2\fv_1 \leftf(\fv, \fv \rightf)
        \right]
        \mathrel +
        \right. \\  
        & \quad
        \left.
        \mathrel +
        \tfrac{1}{24} 
        \left[
        3\:\! \D^2\fv \leftf(\fv, \D\fv\fv \rightf)
        + \D\fv \D\fv \D\fv \fv
        + \D\fv \D^2\fv \leftf( \fv, \fv \rightf)
        + \D^3\fv \leftf(\fv, \fv, \fv \rightf)
        \right]
    \right\rbrace \atyvtt
    \mathrel +
    \ldots 
    \label{eq:taylor3}
\end{align}
The above expression can then be used -- after setting $t=t^j$ -- on the LHS of the condition \eqref{eq:condition}, with the LHS determined by the numerical method $\Phi_{\Dt}$ being investigated. Then \eqref{eq:condition} can yield \emph{recursive relations} for calculating each~$\fv_i$, corresponding to the given numerical method.

\subsection{Demonstration of BEA for the explicit Euler method}\label{sec:ee_example}
First, as an elementary example, let us consider the explicit Euler method for the numerical solution of a general autonomous system \eqref{eq:ode_aut}, written as
\begin{gather}
    \yvjp = \underbrace{\yvj + \Dt\, \fv\leftf( \yvj \rightf)}_{\Phi_{\Dt}\leftf( \yvj \rightf)}.\label{eq:EE}
\end{gather}
Combining \eqref{eq:taylor3} and \eqref{eq:EE} through the condition \eqref{eq:condition} and collecting powers of~$\Dt$ gives the first three of the aforementioned recursive relations, namely,
\begin{align}
    \fv_1 =& - \frac{1}{2} \D \fv \fv, \\
    \fv_2 =& - \frac{1}{2} \D \fv \fv_1 - \frac{1}{2} \D \fv_1 \fv - \frac{1}{6} \D^2\leftf( \fv, \fv \rightf) - \frac{1}{6} \D\fv \D\fv \fv =
    \nonumber
    \\
           =&\ \frac{1}{3} \D \fv \D \fv \fv + \frac{1}{12}\D^2\leftf( \fv, \fv \rightf) , \\
    \fv_3 =& 
        - \frac{1}{2} 
        \left[
            \D\fv_1 \fv_1
            + \D\fv \fv_2
            + \D\fv_2 \fv
        \right]
        -
        \nonumber
        \\ \nonumber
        &- \frac{1}{6} 
        \left[
        2\:\! \D^2\fv \leftf(\fv, \fv_1 \rightf) 
        + \D\fv \D\fv \fv_1
        + \D\fv \D\fv_1\fv
        + \D\fv_1 \D\fv \fv
        + \D^2\fv_1 \leftf(\fv, \fv \rightf)
        \right]
        -
        \\ \nonumber
        &-
        \frac{1}{24} 
        \left[
        3\:\! \D^2\fv \leftf(\fv, \D\fv\fv \rightf)
        + \D\fv \D\fv \D\fv \fv
        + \D\fv \D^2\fv \leftf( \fv, \fv \rightf)
        + \D^3\fv \leftf(\fv, \fv, \fv \rightf)
        \right] =
        \\
        =& 
        -\frac{1}{4} \D\fv \D\fv \D\fv \fv
        -\frac{1}{12} \D\fv \D^2\fv \leftf(\fv, \fv\rightf)
        -\frac{1}{12} \D^2\fv \leftf(\fv, \D\fv \fv\rightf),
\end{align}
which yield the DVF of the explicit Euler method up to and including $\Dt^3$ as
\begin{align}
    \fvt =& \ \fv - \frac{\Dt}{2} \D\fv \fv + \frac{\Dt^2}{12} \D^2\leftf( \fv, \fv \rightf) + \frac{\Dt^2}{3} \D \fv \D \fv \fv - 
    \nonumber \\*
    &-\frac{\Dt^3}{4} \D\fv \D\fv \D\fv \fv
    -\frac{\Dt^3}{12} \D\fv \D^2\fv \leftf(\fv, \fv\rightf)
    -\frac{\Dt^3}{12} \D^2\fv \leftf(\fv, \D\fv \fv\rightf)
    + \ldots 
    \label{eq:EE_modeq}
\end{align}
As can be seen, the DVF obtained equals the original vector field in the zero time-step limit, showing the consistency of the explicit Euler method. Generally, the lowest power of $\Dt$ in the DVF corresponds to the order of the numerical methods \cite{hairer2006geometric}: here, this agrees with the fact that the explicit Euler method is first-order.

It is clearly visible from the above demonstration that the calculation of the DVF can quickly become tedious for even the simplest of all numerical methods as higher-order expansions are considered. A straightforward computer algebra code has been published in \cite{hairer2000asymptotic} for one-dimensional ODEs, which we have generalized for performing BEA on $n$-dimensional systems. This generalized code, which we have used in some of the following calculations, is included in \ref{app:wolfram} for reference.

\subsection{Solution of the truncated distorted equation}
As an illustration of the potential usage of \eqref{eq:EE_modeq}, consider the following linear system corresponding to an undamped 1 degree-of-freedom (DoF) harmonic oscillator (corresponding to a mass-spring system with parameters $m$ and $k$) with natural angular frequency $\omega = \sqrt{k/m}$, position $x(t)$ and velocity $v(t)$:
\begin{gather}
    \underbrace{
    \begin{pmatrix}
        \dot {x} \\
        \dot{v}
    \end{pmatrix}
    }_{\yvdot}
    = 
    \underbrace{
    \begin{pmatrix}
        0 & 1 \\
        -\omega^2 & 0
    \end{pmatrix}
    }_{\mx{A}}
    \underbrace{
    \begin{pmatrix}
        x \\
        v
    \end{pmatrix}
    }_{\yv},
    \label{eq:harmonic_oscillator}
\end{gather}
which can be substituted into \eqref{eq:EE_modeq} as $\fv(\yv) = \mx{A} \yv$, yielding the system of distorted equations
\begin{align}
    \xtdot &= \vt + \frac{\Dt}{2} \omega^2 \xt - \frac{\Dt^2}{3} \omega^2 \vt - \frac{\Dt^3}{4} \omega^4 \xt + \Ord\leftf( \Dt^4 \rightf), \label{eq:harmonic_modeq1}
    \\
    \vtdot &= - \omega^2 \xt + \frac{\Dt}{2} \omega^2 \vt + \frac{\Dt^2}{3} \omega^4 \xt - \frac{\Dt^3}{4} \omega^4 \vt + \Ord\leftf( \Dt^3 \rightf). \label{eq:harmonic_modeq2}
\end{align}

Truncating this after $\Dt^2$ enables us to give a straightforward exact solution to the truncated or approximate distorted equation as
\begin{gather}
    \xt(t) = C_1\;\! \mathrm{e}^{\frac{\Delta t}{2} \omega^2 t} \cos \leftf[ \leftf(\omega -\frac{\Delta t^2}{3} \omega^3 \rightf) t \rightf]+C_2 \;\! \mathrm{e}^{\frac{\Delta t}{2} \omega ^2 t} \sin \leftf[ \leftf(\omega -\frac{\Delta t^2}{3} \omega^3 \rightf) t \rightf],
    \label{eq:harmonic_modeq_sol}
\end{gather}
from which the antidissipative nature of the explicit Euler method is immediately obvious: the amplitude increases over time for $\Dt > 0$. This means that the numerical method introduces additional -- fictional -- energy into the system: a problem for which symplectic numerical methods offer a solution in Hamiltonian systems \cite{devogelaere1956methods,hairer2000asymptotic}. Similarly, a change in the natural angular frequency of the system is also introduced that is on the order of $\Dt^2$.

% Two column figure
\begin{figure}[h]
    \centering
    \includegraphics[width=\textwidth]{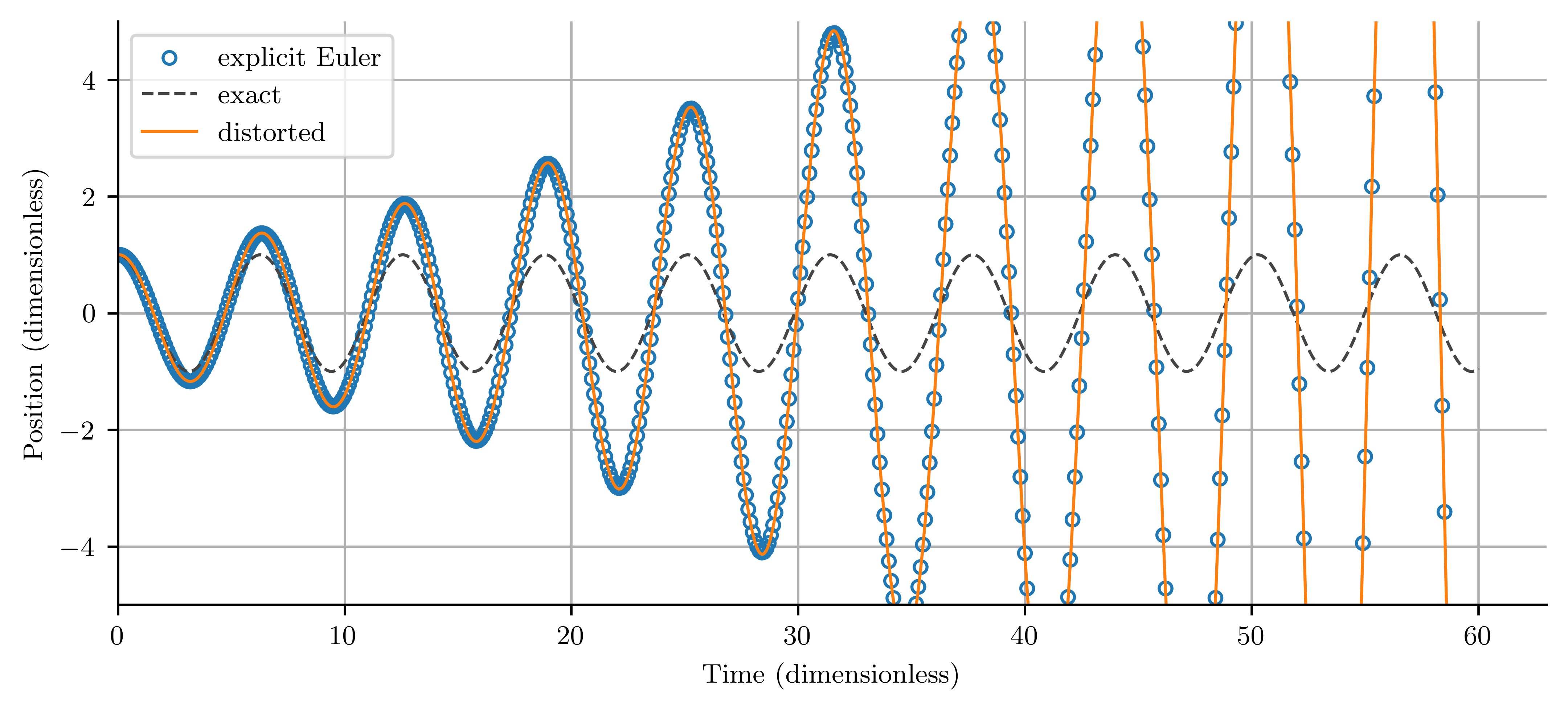}
    \caption{Exact solution, numerical explicit Euler solution and analytical solution of the truncated distorted equation for the position of a mass-spring system as the function of time. Simulation parameters: $\Dt=0.1$, $\omega=0$, $x(0)=1$, $v(0)=0$.}
    \label{fig:EE_solution_time}
\end{figure}

% Single column figure
\begin{figure}[h]
    \centering
    \includegraphics[width=0.5\textwidth]{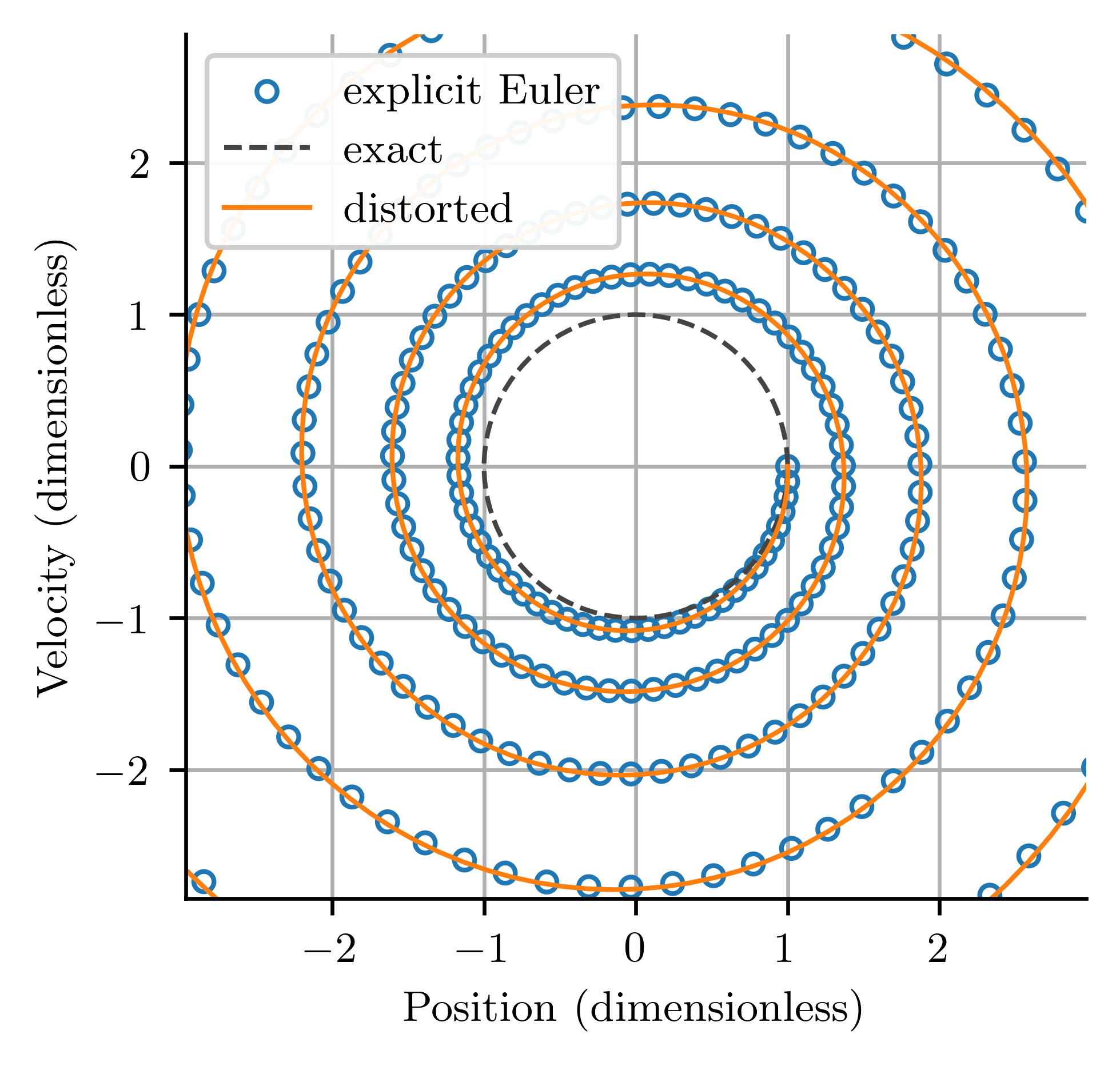}
    \caption{The position and velocity of the results in Fig.~\ref{fig:EE_solution_time} shown in phase space.}
    \label{fig:EE_solution_phase}
\end{figure}

Fig.~\ref{fig:EE_solution_time} shows the comparison between the exact analytical solution and explicit Euler numerical simulation of the original equations~\eqref{eq:harmonic_oscillator}, as well as the solution \eqref{eq:harmonic_modeq_sol} of the truncated distorted equation, for the position as the function of time. Similarly, Fig.~\ref{fig:EE_solution_phase} shows the same results in phase space. It is clearly visible from these graphs that the solution of the truncated distorted equation indeed matches the numerical results closely. Additionally, the antidissipative nature of the distorted system can also be observed.

Additional insight can be gained through the physical interpretation of the distorted equation \eqref{eq:harmonic_modeq1}--\eqref{eq:harmonic_modeq2} by rewriting it into a second-order form with distorted mass, damping and stiffness $\ctilde{m}$, $\ctilde{c}$ and $\ctilde{k}$ as
\begin{gather}
    \ctilde{m} \xtddot + \ctilde{c} \xtdot + \ctilde{k} \xt = 0,
\end{gather}
where
\begin{align}
    \ctilde{m} &= m \left(1+\frac{\Dt^2 k}{3 m} \right) + \Ord\leftf( \Dt^4 \rightf), \\
    \ctilde{c} &= - \Dt \:\! k + \frac{\Dt^3 k^2}{6 m} + \Ord\leftf( \Dt^4 \rightf), \\
    \ctilde{k} &= k \left(1 - \frac{\Dt^2 k}{12 m}\right) + \Ord\leftf( \Dt^4 \rightf),
\end{align}
are the distortions to the original system introduced by the Euler method.

Here we can clearly see that the physical system that the explicit Euler method actually solves is fundamentally different from the original: unlike the original system which had no damping, an additional negative numerical damping is introduced as a numerical artefact. Naturally, the consistency of the method again guarantees that this effect will vanish as the time step approaches zero, but in any practical case where $\Dt > 0$, it persists. This property of the explicit Euler method is well-known; however, using BEA, one can provide an estimate of the amount of numerical damping that is introduced up to a certain order of the time step.

After these motivating examples, we turn to the backward error analysis of the Newmark method.

% A notable and well-known feature of the Newmark method is that for certain parameter combinations, it exhibits no numerical damping. In the following section, we will verify this using BEA, while also deriving further results.

% Here, following [cite moan2006] we assume the existence ???
% - autonóm rendszerre példa?
% - DVF autonóm és nemautonóm részei?

% - explicit Euler, gerjesztetlen rezgő rendszer
% - torzult egyenlet
% - torzult egyenlet egzakt megoldása
% - grafikonok, mint a 2pagerben
% - módosítás?

\section{Distorted equations of the Newmark method}\label{sec:newmark_modified}

\subsection{The Newmark method}
We consider the following system describing the motion of a second-order linear system with $n$ degrees of freedom (represented as $\qv\in\RR^n$) under a non-autonomous external force $\Fv(t)$, written as
\begin{gather}
    \Mm \ddot{\qv}\leftf(t\rightf) + \Cm \dot{\qv}\leftf(t\rightf) + \Km \qv\leftf(t\rightf) = \Fv\leftf(t\rightf), \quad \qv(0)=\qvnull,\, \dot{\qv}(0)=\qdotvnull, \label{eq:structural}
\end{gather}
where $\Mm,\, \Cm,\, \Km\in \RR^{n \times n}$ are the constant mass, damping and stiffness matrices, and $\qvnull$, $\qdotvnull$ are the initial conditions. This represents a linearised but non-autonomous system, which is common in structural analysis, especially in finite element calculations \cite{bathe2014finite,hughes2012finite}, as it is a semidiscrete version of the equation of motion. The Newmark method \cite{newmark1959method, hughes2012finite} is traditionally used for solving \eqref{eq:structural}, along with its various extensions such as the HHT-$\alpha$ \cite{hilber1997improved} and the generalized-$\alpha$ method \cite{chung1993time} (which are usually available in contemporary, state-of-the-art commercial finite element software \cite{ansys2023implicit,abaqus2023implicit}), among many others \cite{wood1980alpha,wilson1968computer,zhou2004design,bathe2005composite,toth2016multifield,malakiyeh2018development,noh2019direct,zakian2021transient,choi2022time,malakiyeh2023explicit}. (More precisely, these are all families of methods with various adjustable parameters, and some of the more general methods contain the less complex ones as special cases.) We only consider systems without a nonlinear internal force in the present paper.

The Newmark method calculates the solution of the generalized coordinate, velocity and acceleration vectors separately ($\qv,\,\vv$ and $\av$, respectively), hence \eqref{eq:structural} is rewritten into a discrete form as
\begin{gather}
    \Mm \avjp + \Cm \vv^{j+1} + \Km \qvjp = \Fvjp, \label{eq:structural_discrete}
\end{gather}
which allows for the expression of $\avjp$ as
\begin{gather}
    \avjp = - \Mminv \left( \Cm \vvjp + \Km \qvjp - \Fvjp\right)
. \label{eq:structural_discrete_acc}
\end{gather}

The two other quantities are calculated for each time step according to the scheme
\begin{align}
    \qvjp &= \qvj + \Dt \:\! \vvj + \frac{\Dt^2}{2}
    \left[ 
        \left(
            1-2\beta
        \right)
        \avj
        +
        2 \beta \avjp
    \right]\label{eq:newmark1}
    \\
    \vvjp &= \vvj + \Dt
    \left[ 
        \left(
            1-\gamma
        \right)
        \avj
        +
        \gamma \avjp
        \right]\label{eq:newmark2}
\end{align}
with $\beta$ and $\gamma$ being the two parameters\footnote{In certain works, the notation $\alpha$ and $\delta$ is also used for the respective parameters.} of the Newmark method given by \eqref{eq:structural_discrete_acc}--\eqref{eq:newmark2}. The choice of these parameters determines the behaviour of the numerical scheme significantly, such as its stability and numerical damping, therefore much research has been focused on the characterization of different parameter combinations. (The Reader is referred to e.g.~\cite{hughes2012finite} for an overview.) However, it must be noted that most of the classical analysis in this topic considers an undamped or specially damped case (i.e.\ proportional or Rayleigh damping, as in \cite{hughes2012finite}) in order to be able to perform a modal transform (mode superposition) and treat the uncoupled 1 DoF timestepping, while we consider the general, $n$-dimensional, linear damping as given in \eqref{eq:structural}, for which direct time integration is usually needed.

Among the special cases of the Newmark family of methods defined by $\beta$ and $\gamma$, of special interest are those where $\gamma=1/2$. These methods are known to produce no numerical damping, a result we will verify using the distorted equations of the Newmark method below. Furthermore, a much less well-known property of the parameter combination $\gamma=1/2$ and $\beta=0$, is it being equivalent (for a Hamiltonian, even nonlinear, variant\footnote{For linear Hamiltonian systems -- such as a second-order linear system without damping and external excitation -- the condition for symplecticity is relaxed to $\gamma=1/2$, $\beta\in \left[0,\frac{1}{2}\right]$}  of \eqref{eq:structural}) \cite{simo1992exact} to the St\"{o}rmer--Verlet method, which is a second-order symplectic numerical scheme \cite{hairer2003geometric}. This result can also be derived from the fact that the Newmark method for non-autonomous systems is, in fact, a variational integrator \cite{kane2000variational}.

\subsection{Deriving the distorted equation}
In order to derive the distorted equation of the Newmark method \eqref{eq:structural_discrete_acc}--\eqref{eq:newmark2} applied to dynamic equation \eqref{eq:structural}, the original system \eqref{eq:structural} needs to be transformed to a suitable form. First we write the equations as a first-order system of ODEs, in a block matrix form as
\begin{gather}
    \begin{pmatrix}
        \dot{\qv} \\ \dot{\vv}
    \end{pmatrix}
    =
    \begin{pmatrix}
        \Nm & \Em \\
        -\Mminv \Km & -\Mminv \Cm
    \end{pmatrix}
    \begin{pmatrix}
        \qv \\ \vv
    \end{pmatrix}
    +
    \begin{pmatrix}
        \Nm \\
        -\Mminv \Fv(t)
    \end{pmatrix}
    \label{eq:structural_nonaut}
\end{gather}
with the newly introduced variable $\vv$
%\equiv 
representing
$\dot{\qv}$, and initial conditions $\qv(0)=\qvnull,\,\vv(0)=\qdotvnull$. However, this is a non-autonomous system due to $\Fv(t)$, and we need an autonomous system to apply the BEA procedure starting with~\eqref{eq:ode_aut}. To achieve this, we can introduce an additional (scalar) degree of freedom $\tau$ as
\begin{gather}
    \underbrace{
    \begin{pmatrix}
        \dot{\tau} \\ \dot{\qv} \\ \dot{\vv}
    \end{pmatrix}
    }_{\dot{\yv}}
    =
    \underbrace{
    \begin{pmatrix}
        1
        \\
        \vv
        \\
        -\Mminv \left( \Cm \vv + \Km \qv - \Fv(\tau) \right)
    \end{pmatrix}
    }_{\fv(\yv)}
    \label{eq:structural_aut}
\end{gather}
with an additional IC as $\tau(0)=0$, and where all the expressions contained in $\dot{\yv}$ are still derivatives of the state variables with respect to $t$ (rather than~$\tau$).

To apply the Newmark method to this extended, nonlinear, autonomous system, the trivial (and exact) timestepping for $\tau$ is needed to supplement \eqref{eq:structural_discrete_acc}--\eqref{eq:newmark2}. After eliminating $\avj$ and $\avjp$ via the substitution of \eqref{eq:structural_discrete_acc}, the full numerical method reads
\begin{align}
    \tau^{j+1} =&\ \tau^j + \Dt, \label{eq:newmark0_2}
    \\
    \qvjp =&\ \qvj + \Dt \:\! \vvj + \frac{\Dt^2}{2}
    \left\lbrace
        \left(
            1-2\beta
        \right)
        \Mminv \left[ \Cm \vvj - \Km \qvj + \Fv\leftf(\tau^j\rightf)\right]
    \right. 
    +
    \nonumber 
    \\
    \label{eq:newmark1_2}
     & \left.
         \mathrel +
        2 \beta 
        \Mminv \left[ \Cm \vvjp - \Km \qvjp + \Fv\leftf(\tau^{j+1}\rightf)\right]
    \right\rbrace,
    \\
    \vvjp =&\ \vvj + \Dt
    \left\lbrace
        \left(
            1-\gamma
        \right)
        \Mminv \left[ \Cm \vvj - \Km \qvj + \Fv\leftf(\tau^j\rightf)\right]
    \right. 
    +
    \nonumber \\ \label{eq:newmark2_2} & \left.
         \mathrel +
        \gamma 
        \Mminv \left[ \Cm \vvjp - \Km \qvjp + \Fv\leftf(\tau^{j+1}\rightf)\right]
        \right\rbrace.
\end{align}
Our goal is now to calculate an asymptotic expansion of the DVF of the method formulated as \eqref{eq:newmark0_2}--\eqref{eq:newmark2_2} when it is applied to the initial value problem~\eqref{eq:structural_aut}. However, the classical approach outlined in Section~\ref{sec:bea_illustration} for constructing the distorted equation has to be extended, since the Newmark method expressed in the form above cannot be written as a function of $\fv$ evaluated at different points in state space -- contrary to the case of most of the traditional numerical methods (including the example in Section~\ref{sec:ee_example}). This is a consequence of the fact that the Newmark method is tailored specifically for problems in structural mechanics, an approach reminiscent of the strategies employed in various structure-preserving numerical methods, where the physical background of the governing equations is also incorporated into the scheme to some extent \cite{shadwick1998exactly, hairer2006geometric, morrison2017structure, fulop2020thermodynamical, takacs2024thermodynamically}.

This obstacle is not significant, though it complicates the calculations somewhat, as the three block components in $\fv(\yv)$ need to be treated separately when dealing with the $\Phi_{\Dt}$ corresponding to the Newmark method. Otherwise, the procedure outlined in Section~\ref{sec:bea_illustration} can be applied, while taking into account that the $2n+1$ DoF system is treated as a block matrix equation, thus the differentiations and series expansions have to be carried out in a scalar ($\tau$) or vectorial ($\qv$, $\vv$) sense for the corresponding block-components.

These calculations were performed using an even more generalized version of the computer algebra code given in \ref{app:wolfram}, which is given in \ref{app:wolfram2}, capable of treating the block matrix equations.

\subsection{Distorted vector field of the Newmark method}
The asymptotic expansion we have obtained for the continuous DVF of the Newmark method \eqref{eq:newmark0_2}--\eqref{eq:newmark2_2} applied to the equation \eqref{eq:structural_aut} is
\begin{gather}
    \fvt = 
    \begin{pmatrix}
        \ctilde{f}_{\tau} \\[1.2ex]
        \fvt_{q} \\[1.2ex]
        \fvt_{v}
    \end{pmatrix}
    =
    \begin{pmatrix}
        1 \\[1.2ex]
        \vv + \Dt\:\! \fvt_{q,1} + \Dt^2 \fvt_{q,2} + \Ord(\Dt^3) \\[1.2ex]
        -\Mminv \left( \Km \qv + \Cm \vv - \Fvftau \right) 
        + \Dt\:\! \fvt_{v,1} + \Dt^2 \fvt_{v,2} + \Ord(\Dt^3)
    \end{pmatrix},
    \label{eq:newmark_mvf}
\end{gather}
where
\begin{align}
    \fvt_{q,1} =&\ \Nm, \label{eq:fvtq1}\\
    \fvt_{q,2} =&\ \eta \Amtqv, \label{eq:fvtq2} \\
    \fvt_{v,1} =&\ \left( \frac{1}{2} - \gamma \right) \Amtqv, \label{eq:fvtv1}\\
    \fvt_{v,2} =&\ \frac{1}{12} \left( 
        \Hm \left(
            \Hm \qv + \Gm \vv - \Mminv \Fvftau
        \right)
        + \Mminv \Fvftaupp
    \right)
    +  \nonumber \\ \label{eq:fvtv2}
               &+ \left( \left( \gamma-\frac{1}{2} \right)^2 + \frac{1}{12} \right) \Gm \Amtqv,
\end{align}
with
\begin{align}
    \eta &= \frac{1}{2} \gamma - \beta - \frac{1}{12}, \\
    \Gm &= \Mminv \Cm, \\
    \Hm &= \Mminv \Km, \\
    \Amtqv &= 
        \left(
        - \Gm \Hm \qv
        +
        \left(
        \Hm - \Gm^2
        \right) \vv
        + \Gm \Mminv \Fvftau
        - \Mminv \Fvftaup
    \right), \label{eq:Amtqv}
\end{align}
and the initial conditions remaining the same as that of the undistorted equation.

From the above expressions, we can quickly verify several properties of the Newmark method. First, it is consistent, as one would expect from any decent numerical method. Furthermore, it is first-order for a generic $\gamma$, as $\fvt_{v,1} \neq \Nv$ per \eqref{eq:fvtv1}. However, a detail that has not been found in the relevant literature is the observation that it is first-order only in one of the variables, and second-order in the other, as $\fvt_{q,1} = \Nv$. It is also immediately apparent from \eqref{eq:fvtv1} that a necessary and sufficient condition of second-order accuracy is $\gamma=1/2$, which is a well-known property of the Newmark method.
%\cite{hughes2012finite}.

As described previously, the presence of numerical damping is dependent on the choice of the parameter $\gamma$. First, consider the undamped case, i.e.\ $\Cm=\Nm$. Then the components of the truncated DVF reduce to
\begin{align}
    \fvt_{q,1} =&\ \Nm, \label{eq:fvtq1_noC}\\
    \fvt_{q,2} =&\ \eta \left( - \Hm \vv + \Mminv \Fvftaup \right), \label{eq:fvtq2_noC} \\
    \fvt_{v,1} =&\ \left( \gamma - \frac{1}{2} \right) \left( - \Hm \vv + \Mminv \Fvftaup \right), \label{eq:fvtv1_noC}\\
    \fvt_{v,2} =&\ \frac{1}{12} \left( 
        \Hm \left(
            \Hm \qv - \Mminv \Fvftau
        \right)
        + \Mminv \Fvftaupp
    \right)
    \label{eq:fvtv2_noC}
\end{align}
where the only numerical damping remaining is the term containing $\vv$ in \eqref{eq:fvtv1_noC}. By setting $\gamma=1/2$ this whole component vanishes, thus numerical damping is indeed eliminated, as only $\qv$- and $\Fv$-dependent terms remain in $\fvt_{v,2}$. This agrees with the results described in \cite{simo1992exact} for conservative linear Hamiltonian systems, as the absence of numerical damping does not depend on $\beta$ for such systems.

On the other hand, this is not the case for $\Cm \neq \Nm$. As it is visible from \eqref{eq:fvtq1}--\eqref{eq:fvtv2}, though setting $\gamma=1/2$ does eliminate $\vv$-dependent terms from $\fvt_{v,1}$, it does not eliminate such terms from $\fvt_{v,2}$, which correspond to additional numerical damping on top of the physical damping created by $\Cm$. This means that, for damped systems, setting $\gamma=1/2$ does not fully eliminate numerical damping. 

For the numerical verification and illustration of these results, the Reader is referred to Section~\ref{sec:num_verification} and Section \ref{sec:num_demo1}, respectively. 

Before the illustrations, however, we will use the present results to relate the distorted vector field to the original, physical, interpretation of the system equations.

\subsection{Distorted second-order equation corresponding to the Newmark method}
Equation \eqref{eq:newmark_mvf} can be transformed back into the original second-order form, which can be written as an equation of motion with distorted coefficients and forcing, namely,
\begin{gather}
    \Mm \ddot{\qv}(t) + \Ctm \dot{\qv}(t) + \Ktm \qv(t) = \Ftvft, \label{eq:newmark_mod_secondorder}
\end{gather}
where the distorted matrices and forcing are
\begin{align}
    \Ctm =& \ 
    \Cm + 
    \Bmfhg
    \left(
         \Km - \Cm \Mminv \Cm
    \right)
    + \Dt^2 \left(\eta - \frac{1}{12}\right) \Km \Mminv \Cm
    + \Ord\leftf( \Dt^3 \rightf),
    \label{eq:Ctm}
    \\*
    \Ktm =& \
    \Km -
    \Bmfhg \Cm \Mminv \Km
    + \Dt^2 \left(\eta - \frac{1}{12}\right) \Km \Mminv \Km
    + \Ord\leftf( \Dt^3 \rightf),
    \label{eq:Ktm}
    \\*
    \Ftvft =& \
    \Fvft 
    - \Bmfhg \left( \Cm \Mminv \Fvft - \Fvftp \right)
    + \Dt^2 \left(\eta - \frac{1}{12}\right)\left( \Km \Mminv \Fvft - \Fvftpp \right),
    \label{eq:Ftvft}
\end{align}
up to second order, with
\begin{align}
    \eta &= \frac{1}{2} \gamma - \beta - \frac{1}{12}, \\
    \Bmfhg &= 
            \Dt \left( \gamma - \frac{1}{2} \right) \Em 
            -
            \Dt^2 \left( \left( \gamma - \frac{1}{2} \right)^2 + \frac{1}{12} \right) \Cm \Mminv
        .
\end{align}

One can observe that due to the nature of the rewriting \eqref{eq:newmark_mod_secondorder}, $\dot{\qv}(t)$ there does not equal the $\vv(t)$ in \eqref{eq:newmark_mvf}. Consequently, the initial condition for the first-order derivative changes compared to the original velocity initial condition, according to the second component in \eqref{eq:newmark_mvf}. Accordingly, the initial conditions of \eqref{eq:newmark_mod_secondorder} become
\begin{align}
    \qv(0) &= \qvnull, \\
    \dot{\qv}(0) &= \qdotvnull + \Dt^2 \eta \Am\leftf( \tau{=}0, \qv{=}\qvnull, \vv{=}\qdotvnull \rightf).
\end{align}

The distorted matrices enable us to characterize the behaviour of the numerical method in terms of parameters of the distorted system. For $\gamma=1/2$, the distorted damping matrix reduces to
\begin{gather}
    \eval{\Ctm}_{\gamma=1/2} =
    \Cm + \frac{1}{12} \Dt^2 \left[ 
        \left(
            1 - 12 \beta
        \right)
        \Km \Mminv \Cm
        - \Cm \Mminv \left( \Km -\Cm \Mminv \Cm \right)
    \right]
    + \Ord\leftf( \Dt^3 \rightf).
\end{gather}
Indeed, in line with the analysis above, for $\Cm=\mx{0}$ this expression does equal zero, indicating the absence of numerical damping -- for at least second order --, while for the damped case of $\Cm \neq \mx{0}$, there is some (second-order) numerical damping remaining. This stems from the distortion of the original damping matrix by the Newmark method in the distorted equation, which cannot be recognized using any analysis of the undamped case. Additionally, it should be noted that this result does not contradict the previous statements about the symplecticity of the Newmark method -- as symplecticity can only be interpreted for conservative Hamiltonian systems, while a damped system is clearly not conservative.

\subsection{Numerical verification}\label{sec:num_verification}
We can use numerical tools for verifying the correctness of the derived distorted vector field and the distorted system parameters. According to \eqref{eq:newmark_mvf} and \eqref{eq:newmark_mod_secondorder}, the distorted vector field and second-order equation are accurate up to and including $\Ord\leftf( \Dt^2 \rightf)$, thus the deviation between the numerical results from a Newmark numerical simulation and the solution of \eqref{eq:newmark_mvf} or \eqref{eq:newmark_mod_secondorder} should both be of~$\Ord\leftf( \Dt^3 \rightf)$.

For showing this, we use a 3 DoF-system with (randomly generated) parameters
\begin{gather}
    \Mm = 
    \begin{pmatrix}
        4.6965 & 1.4187 & 1.6038 \\
        1.4187 & 4.7195 & 1.5540 \\
        1.6038 & 1.5540 & 4.4809 \\
    \end{pmatrix}
    ,\quad
    \Km = 
    \begin{pmatrix}
        4.5316 & 1.6906 & 1.6784 \\
        1.6906 & 4.7245 & 1.4670 \\
        1.6784 & 1.4670 & 4.3618 \\
    \end{pmatrix},
    \\
    \Cm = 
    \begin{pmatrix}
        0.033921 & 0.003909 & 0.007335 \\
        0.003909 & 0.030597 & 0.002903 \\
        0.007335 & 0.002903 & 0.031755 \\
    \end{pmatrix},\label{eq:C_param}
    \\
    \Fvft = 
    \begin{pmatrix}
        -0.040790 \cdot \sin(0.2457\cdot t) \\ 
        -0.006630 \cdot \sin(0.2587\cdot t) \\ 
        -0.006914 \cdot \sin(0.3262\cdot t) \\
    \end{pmatrix},
    \label{eq:sys11}
\end{gather}    
and initial conditions
\begin{gather}
    \qvnull = 
    \begin{pmatrix}
        0.1 \\ 0 \\ 0 \\ 
    \end{pmatrix}
    , \quad 
    \qdotvnull =
    \begin{pmatrix}
        0 \\ 0 \\ 0 \\ 
    \end{pmatrix}.
    \label{eq:sys12}
\end{gather}
The simulations were ran using progressively finer time steps $\Dt$, and were all evaluated at $t = 0.4$. For the Newmark simulations, the parameters were $\gamma=0.55$ and $\beta=0.28$ for achieving as general a comparison as possible. For solving \eqref{eq:newmark_mvf} and \eqref{eq:newmark_mod_secondorder} numerically, a fixed-step, fourth-order Runge--Kutta (RK4) solver has been used with time step $\Dt/100$.

\begin{figure}[h]
    \centering
    \includegraphics[width=0.495\textwidth]{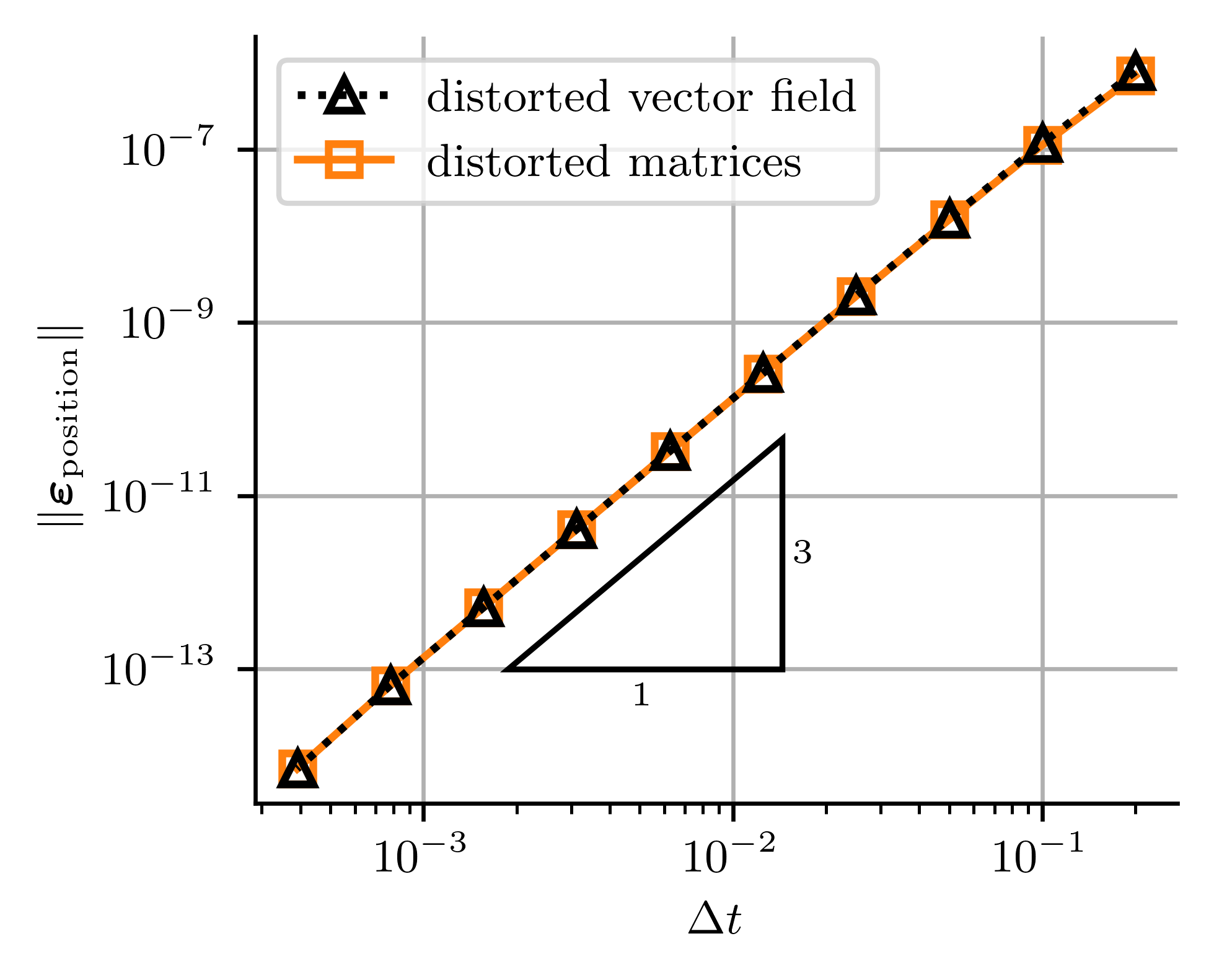}
    \hfill
    \includegraphics[width=0.495\textwidth]{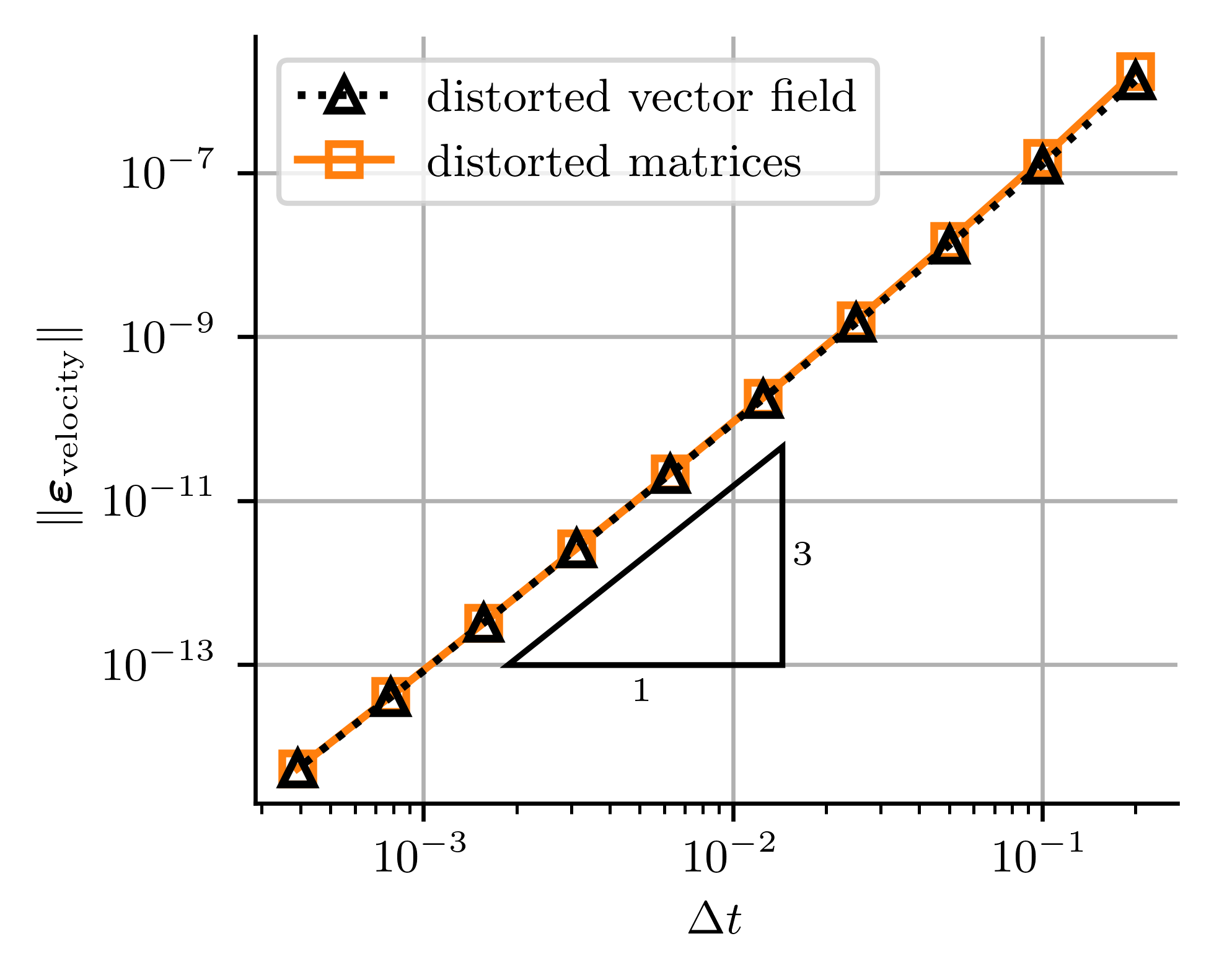}
    \caption{Convergence of the norm of deviation between the Newmark simulation results and the distorted vector field or the second-order equation, for the position (left) and velocity (right) solutions.}
    \label{fig:DVF_mat_error}
\end{figure}

Fig.~\ref{fig:DVF_mat_error} shows for the position and velocity solutions that both the derived distorted vector field and the second-order equation parameters are correct, as the deviations achieve the desired convergence.

% Two column figure

\subsection{Numerical demonstration}\label{sec:num_demo1}
\begin{figure}[h]
    \centering
    \includegraphics[width=\textwidth]{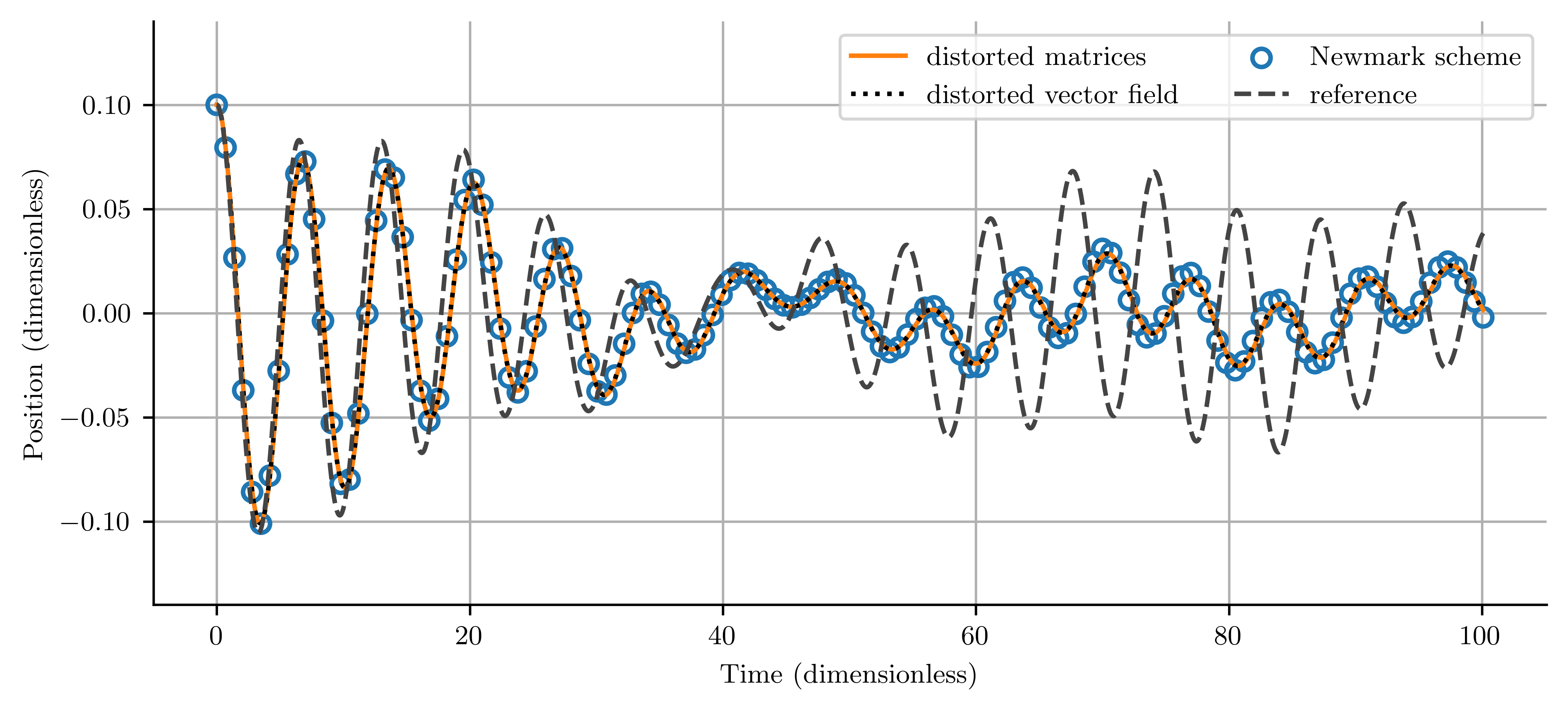}
    \caption{Solution for the first component of the position vector $\qv$ as the function of time, using the Newmark method, compared to solutions of the distorted vector field \eqref{eq:newmark_mvf} and distorted second-order equation \eqref{eq:newmark_mod_secondorder}, as well as a reference solution of the original system. (System parameters per \eqref{eq:sys11}--\eqref{eq:sys12}, $\Dt=0.7$.)}
    \label{fig:DVF_mat_time}
\end{figure}

The agreement between the numerical results from a Newmark simulation and the two distorted equations \eqref{eq:newmark_mvf} and \eqref{eq:newmark_mod_secondorder} can also be visualized as the function of time, as shown in Fig.~\ref{fig:DVF_mat_time} along with a numerical reference solution\footnote{As a numerical \emph{reference solution}, in this paper we refer to a fourth-order Runge--Kutta solution that is calculated with at least one hundredth of the timestep used for obtaining the other numerical results in a certain comparison. Consequently, such a reference solution serves as an adequate substitute for an exact solution in these comparisons.} of the original system. A clear difference can be observed between the reference solution and the Newmark solution, partly due to the relatively large time step used (frequency mismatch), and partly due to the numerical damping introduced (amplitude mismatch).

Meanwhile, the solutions of both continuous, distorted equations match the results of the Newmark simulation closely. This again shows that they capture the behaviour of the discrete numerical method in continuous time accurately, both in a quantitative and qualitative sense.

% Two column figure

\section{Applications of the distorted equations for improving the Newmark method}\label{sec:newmark_improvement}
In the above, we have established that a numerical method applied to a system of ODEs modify the system in the sense that the numerical results are an exact solution of a different, underlying distorted system of ODEs, and the distorted equation can be approximated by an asymptotic series that can be calculated according to techniques employed in BEA. We have derived the distorted system of ODEs corresponding to the Newmark method applied to the second-order, linear structural equation with non-autonomous forcing, and have shown that this result can be interpreted as an underlying second-order system with distorted damping, stiffness and forcing.

While this novel result already contains useful information for analysis of the Newmark method, we will show in the following section that the distorted equations can also be used in a constructive way to achieve better performance or accuracy in simulations using the Newmark method -- without
%directly
modifying the numerical scheme itself, or changing the size of the time step, in an approach we call backward error analysis-based \emph{\compensation}.

The central idea of 
%the backward error analysis-based
this 
\compensation
%(or simply \emph{compensation} for short, from here on) 
is that by knowing the distortions introduced by a numerical method to the differential equation, the physical parameters of the system being investigated can be \compensated accordingly, in order to obtain more accurate results in a certain desired sense. This is performed in such a way that the distortions introduced by the numerical method are cancelled out (at least up to a chosen order of $\Delta t$) by the \compensation, thus the results from the numerical methods are actually closer to the exact results than if the simulation were performed with the original system parameters.

This approach bears some similarity to the exactly conservative integrators of Shadwick et al \cite{shadwick1998exactly} and the distorted Hamiltonian of the Takahashi--Imada integrator \cite{takahashi1984monte}, as well as the mYBABY method \cite{shang2020structurepreserving}, though there is an important distinction: in those cases, the system equation or the numerical method itself is modified in order to achieve better numerical results, while in the case of the \compensation introduced here, the modification of the system parameters (including the forcing function in certain cases) is sufficient. This distinction is especially significant if one considers existing implementations of numerical simulation methods: compensation can be applied to the input parameters of the time integration, thus superior results could be achieved using existing, potentially proprietary software.

In what follows, we introduce two \compensations of the Newmark method: one deals with the elimination of the numerical damping, and the other achieves a fourth-order accuracy using the originally second-order accurate Newmark method.

\subsection{Eliminating numerical damping from the Newmark method using a \compensated damping matrix}
The expression \eqref{eq:Ctm} gives a second-order estimate of the distorted damping matrix. As discussed above, this shows that numerical damping is introduced by the Newmark method into the simulation results, the extent of which also depends on the value of the Newmark parameter $\gamma$. For \compensating this aspect of the numerical simulation, a \compensated damping matrix can be introduced, that guarantees the absence of numerical damping, up to a certain order of the time step.

According to the above, the condition to be fulfilled for the \compensated damping matrix $\Chm$ is
\begin{gather}
    \Cm = \ 
    \Chm + 
    \Bmfhg
    \left(
         \Km - \Chm \Mminv \Chm
    \right)
    + \Dt^2 \left(\eta - \frac{1}{12}\right) \Km \Mminv \Chm
    + \Ord\leftf( \Dt^3 \rightf),\label{eq:Chm_condition}
\end{gather}
where we look for $\Chm$ in the following form:
\begin{gather}
    \Chm = \Cm + \Dt \Chm_1 + \Dt^2 \Chm_2.\label{eq:Chm_definition}
\end{gather}
The zeroth-order term is set to $\Cm$ as we would like to maintain the consistency of the \compensated method, and the remaining terms are introduced up to second order of the time step, in line with the original expression. Inserting \eqref{eq:Chm_definition} into \eqref{eq:Chm_condition}, expanding and collecting terms by order of time step yields the equations
\begin{align}
    \Dt^1: \: 0=&\ \Chm_1 + \left( \gamma - \frac{1}{2} \right) \left(\Km - \Cm \Mminv \Cm \right),
    \label{eq:Chm1_int}
    \\
    \Dt^2: \: 0=&\ \Chm_2 - \left( \gamma - \frac{1}{2} \right)^2 \left( \left(\Cm \Mminv \Cm - \Km \right)\Mminv \Cm  + \Cm \Mminv \left( \Cm \Mminv \Cm - \Km \right) \right) - 
    \nonumber
    \\
    & -\left(\left( \gamma - \frac{1}{2} \right)^2 + \frac{1}{12} \right) \Cm \Mminv \left( \Km - \Cm \Mminv \Cm \right) + \left( \frac{1}{2}\gamma - \beta - \frac{1}{6} \right) \Km \Mminv \Cm.
\label{eq:Chm2_int}
\end{align}
Solving \eqref{eq:Chm1_int} and \eqref{eq:Chm2_int} gives the value of the two \compensating terms as
\begin{align}
    \Chm_1 =& \ \left( \gamma - \frac{1}{2} \right) \left(\Cm \Mminv \Cm  - \Km \right),\label{eq:Chm1}
    \\
    \Chm_2 =& \ \left(\left( \gamma - \frac{1}{2} \right)^2 - \frac{1}{12} \right) \Cm \Mminv \Cm \Mminv \Cm \mathrel - \nonumber \\
            & \ - \left( \gamma^2 - \frac{1}{2} \gamma - \beta + \frac{1}{12} \right) \Km \Mminv \Cm + \frac{1}{12} \Cm \Mminv \Km.\label{eq:Chm2}
\end{align}

\subsubsection{Numerical demonstration (undamped case)}\label{sec:num_demo2_compC_noC}
To demonstrate the numerical damping compensation, the system described in Section~\ref{sec:num_verification} is used, without any excitation (i.e.~$\Fvft \equiv \mx{0}$), and first with no damping (i.e.~$\Cm=\mx{0}$), to make the effects of numerical damping more apparent.

Using \eqref{eq:Chm_definition}, \eqref{eq:Chm1} and \eqref{eq:Chm2}, the compensated damping matrix $\Chm$ for the Newmark method can be calculated, with the same $\Dt$, $\gamma=0.55$ and $\beta=0.28$ as in Section~\ref{sec:num_verification}. The compensated damping matrix eliminates the numerical damping introduced by the Newmark method. In other words, more accurate simulation of an undamped system can be achieved by simulating a fictional compensated system that has an additional (in this case, negative) compensating damping introduced.

For comparison, the performance of both the fourth-order Runge--Kutta method and the generalized-$\alpha$ method are given in the following examples. For the latter method, $\rho_{\infty}=0.9$ is used (as a typical value), which also yields $\gamma=0.55$ and $\beta=0.28$ using the formulas for the optimal $\alpha_f$, $\alpha_m$, $\gamma$ and $\beta$ values \cite{chung1993time}, making the comparison as fair as possible.

\begin{figure}[h]
    \centering
    \includegraphics[width=\textwidth]{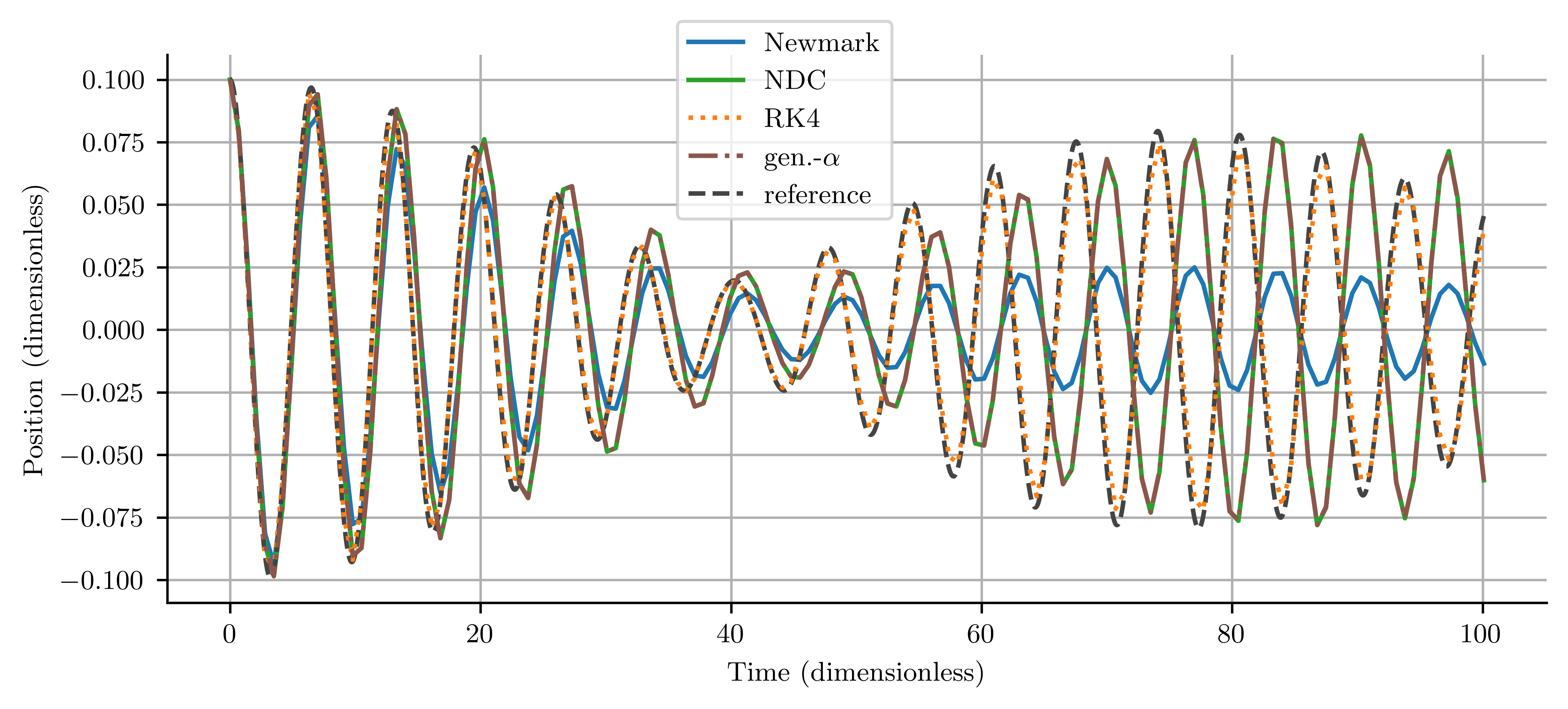}
    \caption{Solution for the first component of the position vector $\qv$ as the function of time, from simulating the original, undamped system using the Newmark method, compared to the Newmark damping compensated (NDC) system also simulated with the Newmark method, against a Runge--Kutta (RK4), generalized-$\alpha$ and a reference solution. (Time step for the former four was identically $\Dt=0.7$.)}
    \label{fig:compC_position_time_t100_noC}
\end{figure}

The validity of the approach is demonstrated using numerical results. Fig.~\ref{fig:compC_position_time_t100_noC} shows the simulated position as the function of time: it is clearly visible that the numerical damping is indeed eliminated, while the frequency mismatch introduced by the Newmark method still remains. This highlights the possibility to use backward error analysis-based compensation in a selective manner to only compensate certain aspects of a numerical method, while leaving others untouched. Remarkably, the damping-compensated Newmark method gives similar position results to the generalized-$\alpha$ method.

% Two column figure

\begin{figure}[h]
    \centering
    \includegraphics[width=0.495\textwidth]{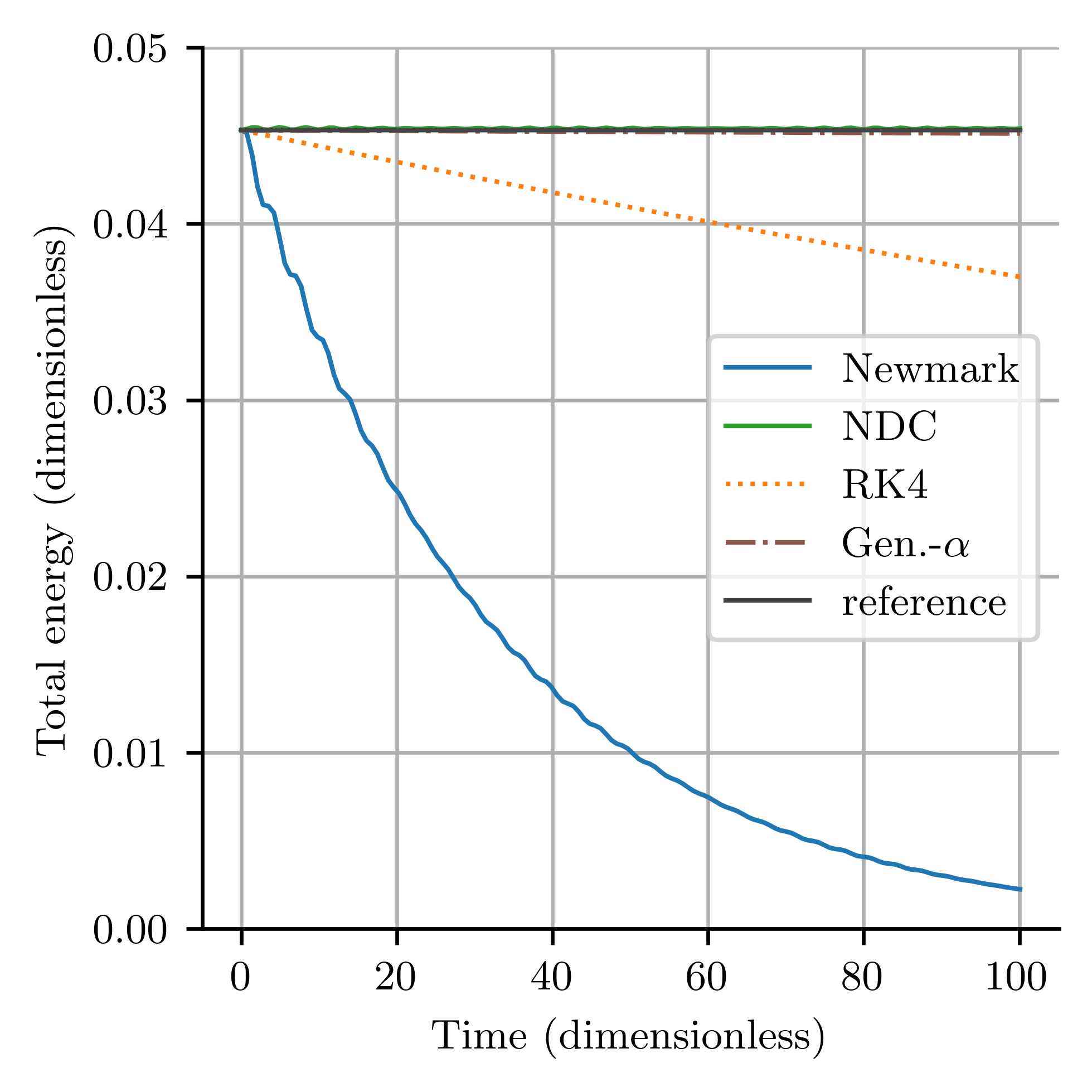}
    \includegraphics[width=0.495\textwidth]{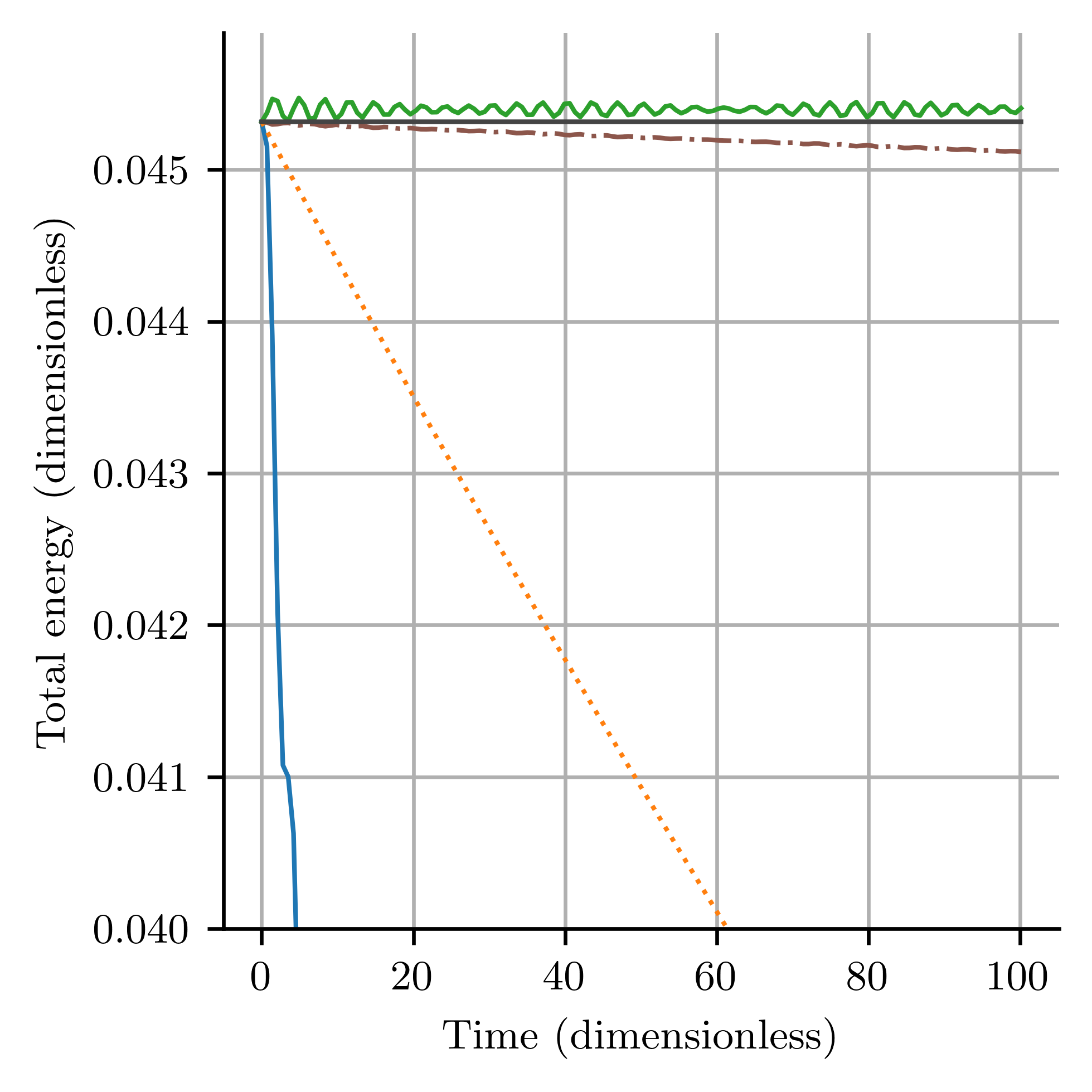}
    \caption{Total energy as the function of time, from simulating the original undamped system using the Newmark method, compared to the Newmark damping compensated (NDC) system also simulated with the Newmark method, against a Runge--Kutta (RK4), generalized-$\alpha$ and a reference solution. (Time step for the former four was identically $\Dt=0.7$.)}
    \label{fig:compC_E_tot_time_t200_noC}
\end{figure}

The same results can be used to demonstrate the elimination of numerical damping from an energy perspective, as shown in Fig.~\ref{fig:compC_E_tot_time_t200_noC}. Here the almost-constant nature of the total energy in the compensated system is clearly visible, while both the Runge--Kutta and the uncompensated Newmark scheme simulations show exponential numerical dissipation, with the generalized-$\alpha$ method also exhibiting a slight numerical dissipation. The small oscillations in the total energy of the compensated system are due to the fact that the numerical damping compensation is only second order and not exact, and are visible due to the large time step: with decreasing the time step, the oscillations in the total energy also decrease.

% Two column figure

\subsubsection{Numerical demonstration (damped case)}\label{sec:num_demo2_compC_withC}
For showing the validity of the numerical damping compensation, the same approach is also applied to a damped system. This might be even more relevant for the application of this approach, as the numerical damping of the Newmark scheme is sometimes beneficial in the simulation of an undamped system, but if physical damping is already present, numerical damping is significantly less relevant, and its effects can be confused with the effects of the physical damping.

\begin{figure}[h]
    \centering
    \includegraphics[width=\textwidth]{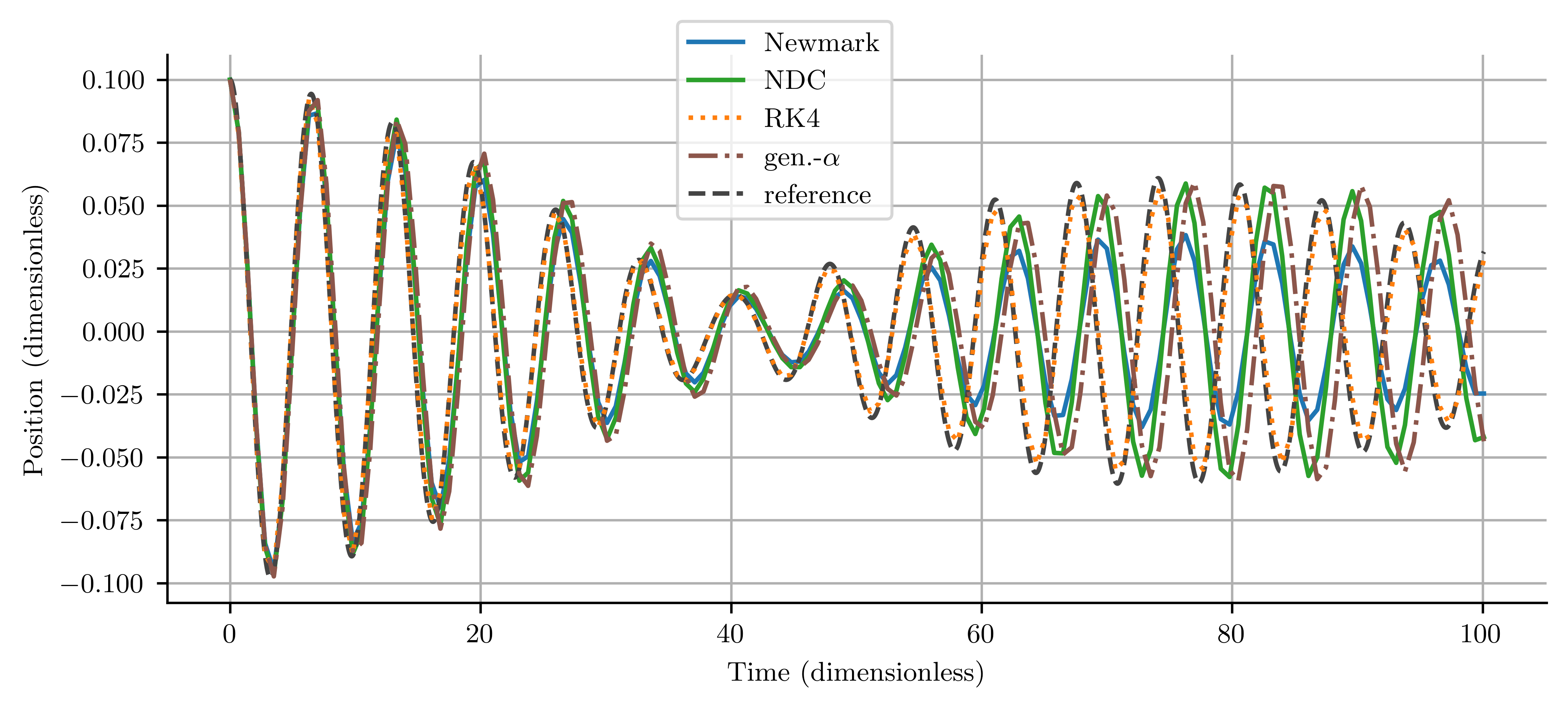}
    \caption{Solution for the first component of the position vector as the function of time for the damped 3 DoF system with no excitation. Newmark damping compensation (NDC) against uncompensated Newmark, RK4, generalized-$\alpha$ and reference solutions.}
    %\caption{Solution for the first component of the position vector $\qv$ as the function of time, from simulating the original, damped system using the Newmark method, compared to the Newmark damping compensated (NDC) system also simulated with the Newmark method, against a Runge--Kutta (RK4) and a reference solution. (Time step for the former three was identically $\Dt=0.7$.)}
    \label{fig:compC_position_time_t100}
\end{figure}

\begin{figure}[h]
    \centering
    \includegraphics[width=0.5\textwidth]{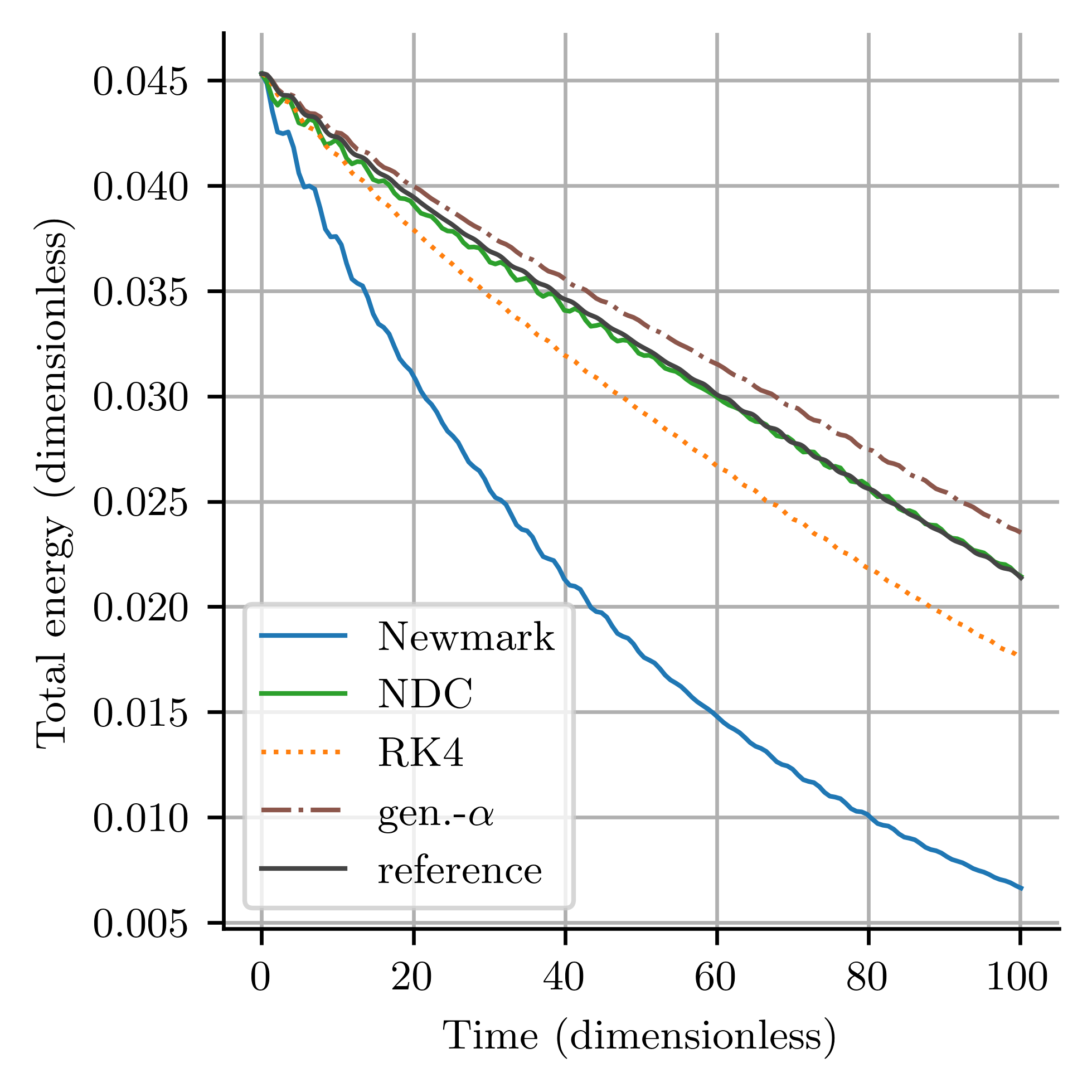}
    %\caption{Total energy as the function of time, from simulating the original, damped system using the Newmark method, compared to the Newmark damping compensated (NDC) system also simulated with the Newmark method, against a Runge--Kutta (RK4) and a reference solution. (Time step for the former three was identically $\Dt=0.7$.)}
    \caption{Total energy as the function of time for the damped 3 DoF system with no excitation. Newmark damping compensation (NDC) against uncompensated Newmark, RK4, generalized-$\alpha$ and reference solutions.}
    \label{fig:compC_E_tot_time_t200}
\end{figure}

In this demonstration, the same setup has been used as in Section~\ref{sec:num_demo2_compC_noC}, with the only difference that now the physical damping is nonzero, having the value given in \eqref{eq:C_param}. Figs.~\ref{fig:compC_position_time_t100} and~\ref{fig:compC_E_tot_time_t200} show the position and total energy solution as the function of time, respectively. The damping compensation again eliminates almost all numerical damping, while respecting the physical damping of the system. Meanwhile, both the uncompensated Newmark and the fourth-order Runge--Kutta schemes show clear signs of numerical dissipation. Additionally, the generalized-$\alpha$ method shows a different phase compared to the Newmark and compensated Newmark methods, while introducing an overall negative, unphysical, numerical dissipation.

% Two column figure

% Single column figure

Based on the above, we can conclude that the numerical damping compensation of the Newmark method presented here achieves the goal of more accurately representing the original system under the Newmark scheme. Moreover, this is achieved while using the same amount of resources for the simulation itself, after the one-time calculation of the compensated damping matrix $\Chm$.

In the following subsection, we further extend this approach.

\subsection{Fourth-order accuracy with the Newmark method using \compensated parameters}
We have already stated previously that the lowest nonzero-order term in a DVF corresponds to the order of the numerical method itself. This property of distorted equations lends itself to the idea that by cancelling out the lower-order time step terms during \compensation, the accuracy of a numerical method could be increased. In what follows, we will show the process of obtaining \compensation parameters to make the Newmark method fourth-order.

It is known that, for $\gamma=1/2$, the Newmark-method is symmetric \cite{bathe2014finite,hughes2012finite}. A property of symmetric numerical schemes is that they only have even-order time-step terms in the asymptotic expansions of their distorted equations \cite{hairer2006geometric}. Thus, if $\gamma=1/2$ is set for the Newmark method, and the $\Dt^2$ term in the DVF is cancelled using \compensation, then the lowest-order remaining term will be of $\Dt^4$, and thus the \compensated method will be fourth-order accurate.

The ansatz for the \compensated matrices and forcing is similar to the previous construction, but without first-order terms, in line with the absent odd-order terms in the DVF: 
\begin{align}
    \Chm &= \Cm + \Dt^2 \Chm_2, \label{eq:rect4Chm}\\
    \Khm &= \Km + \Dt^2 \Khm_2, \label{eq:rect4Khm}\\
    \Fhvft &= \Fvft + \Dt^2 \Fhvtwoft, \label{eq:rect4Fhvft}
\end{align}
where the function $\Fhvtwoft$ might be dependent on derivatives of $\Fvft$.

To determine the \compensating terms, we substitute $\gamma=1/2$ into the DVF \eqref{eq:newmark_mvf}. Observe that the component $\fvt_{q}$ becomes
\begin{gather}
    \eval{\fvt_{q}}_{\gamma=\frac{1}{2}} = 
    \vv + 
    \Dt^2 \left( \frac{1}{6} - \beta \right) \Amtqv + \Ord\leftf( \Dt^4 \rightf),
\end{gather}
and as -- due to \eqref{eq:Amtqv} -- only the zeroth-order terms in \eqref{eq:rect4Chm}--\eqref{eq:rect4Fhvft} can influence the second-order term in this expression, we find that the constraint $\beta=1/6$ is also necessary to achieve fourth-order accuracy, besides the already imposed constraint on $\gamma$ (i.e., $\gamma=1/2$). With these constraints, the DVF \eqref{eq:newmark_mvf} reduces to
\begin{gather}
    \eval{
        \begin{pmatrix}
            \ctilde{f}_{\tau} \\
            \fvt_{q} \\
            \fvt_{v}
        \end{pmatrix}
    }_{\gamma=\frac{1}{2}, \beta=\frac{1}{6}}
    =
    \begin{pmatrix}
        1 \\
        v + \Ord\leftf( \Dt^4 \rightf) \\
        -\Mminv \left( \Km \qv - \Cm \vv + \Fvftau \right) 
        + \eval{\Dt^2 \fvt_{v,2}}_{\gamma=\frac{1}{2}, \beta=\frac{1}{6}}
        + \Ord(\Dt^4) \\
    \end{pmatrix},
\end{gather}
thus only $\fvt_{v}$ remains to be \compensated. Eliminating the dummy variable $\tau$, the condition for this is
\begin{align}
    -\Mminv \left( \Km \qv - \Cm \vv \right. & \left. + \Fvft \right) = -\Mminv \left( \Khm \qv - \Chm \vv + \Fhvft \right) + \nonumber \\*
    & \qquad + \frac{1}{12} \Dt^2 \left[ 
        \left(\Mminv\Khm - \Mminv\Chm\Mminv\Chm \right)\Mminv \Khm \qv
        \mathrel + \nonumber \right. \\
    & \qquad + \left(
    \Mminv \Chm \left(\Mminv\Khm - \Mminv\Chm\Mminv\Chm \right) + \Mminv \Khm \Mminv \Chm 
        \right) \vv
        \mathrel - \nonumber \\
    & \qquad - \left(\Mminv\Khm - \Mminv\Chm\Mminv\Chm \right) \Mminv \Fhvft 
        \mathrel - \nonumber \\
        & \qquad \left.  \mathrel - \Mminv \Chm \Mminv \Fhvftp + \Mminv \Fhvftpp
    \right]
    + \Ord(\Dt^4)
    ,
\end{align}
which, after substitution of the ansatz \eqref{eq:rect4Chm}--\eqref{eq:rect4Fhvft}, and subsequent expansion and collection of terms, yields as its solution the \compensating terms
\begin{align}
    \Chm_2 &= \frac{1}{12} \left( 
        \Cm \Mminv \Km + \Km \Mminv \Cm - \Cm \Mminv \Cm \Mminv \Cm
    \right),
    \\
    \Khm_2 &= \frac{1}{12} \left(
        \Km \Mminv \Km - \Cm \Mminv \Cm \Mminv \Km
    \right),
    \\
    \Fhvtwoft &= \frac{1}{12} \left( 
        \Cm \Mminv \left(
            \Cm \Mminv \Fvft - \Fvftp
        \right)
        - 
        \Km \Mminv \Fvft + \Fvftpp
    \right).
\end{align}

Thus, the \compensated matrices and excitation for a fourth-order accurate calculation are
\begin{align}
    \Chm &= \Cm + \frac{1}{12} \Dt^2 \left( 
        \Cm \Mminv \Km + \Km \Mminv \Cm - \Cm \Mminv \Cm \Mminv \Cm
    \right),\label{eq:Chm}
    \\
    \Khm &= \Km + \frac{1}{12} \Dt ^2 \left(
        \Km \Mminv \Km - \Cm \Mminv \Cm \Mminv \Km
    \right),\label{eq:Khm}
    \\
    \Fhvft &= \Fvft + \frac{1}{12} \Dt^2 \left( 
        \Cm \Mminv \left(
            \Cm \Mminv \Fvft - \Fvftp
        \right)
        - 
        \Km \Mminv \Fvft + \Fvftpp
    \right).
    \label{eq:Fhvft}
\end{align}

\subsubsection{Numerical demonstration (1 DoF convergence)}\label{sec:convergence1dof}

As the first demonstration of the fourth-order compensation system, we investigate the convergence of a compensated system simulated using the Newmark scheme. For this, we use a 1 DoF system with the same parameters as in \cite{hulbert1987error}, with an added harmonic excitation, for which a closed-form exact solution exists and can be used as an accurate reference. The parameters of the system are
\begin{gather}
    \Mm = m,\,\quad \Cm = 2 \xi \omega,\,\quad \Km=m \omega^2, \\
    \Fvft = 0.8 \cos( 10 \omega ),
\end{gather}
where
\begin{gather}
    m = 1,\,\quad \xi=0.02,\,\quad \omega = 2 \pi,
\end{gather}
are the mass, damping parameter and angular frequency, respectively. The initial conditions are
\begin{gather}
    \qvnull = 
    1
    , \, 
    \qdotvnull =
    1.
\end{gather}

\begin{figure}[h]
    \centering
    \includegraphics[width=0.495\textwidth]{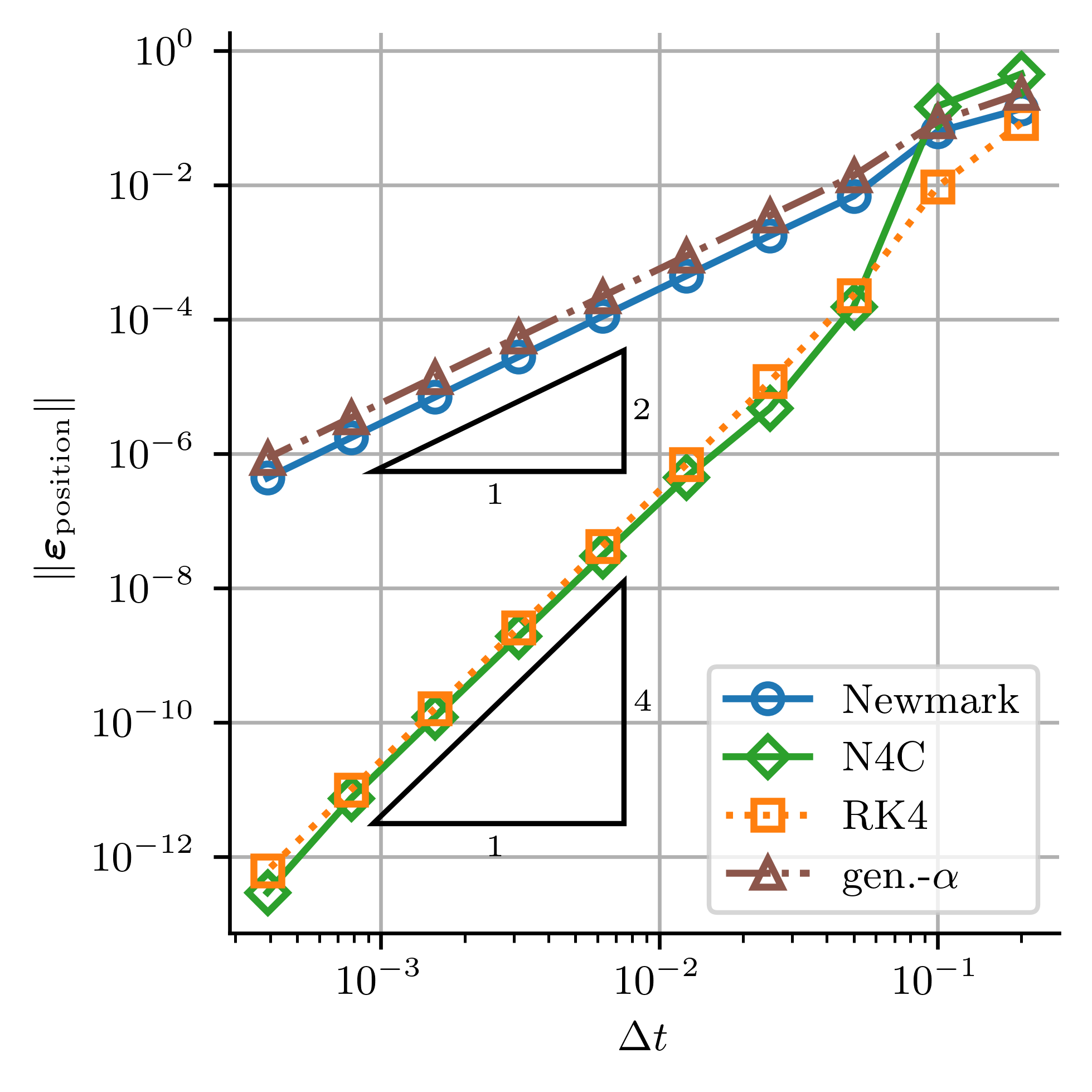}
    \includegraphics[width=0.495\textwidth]{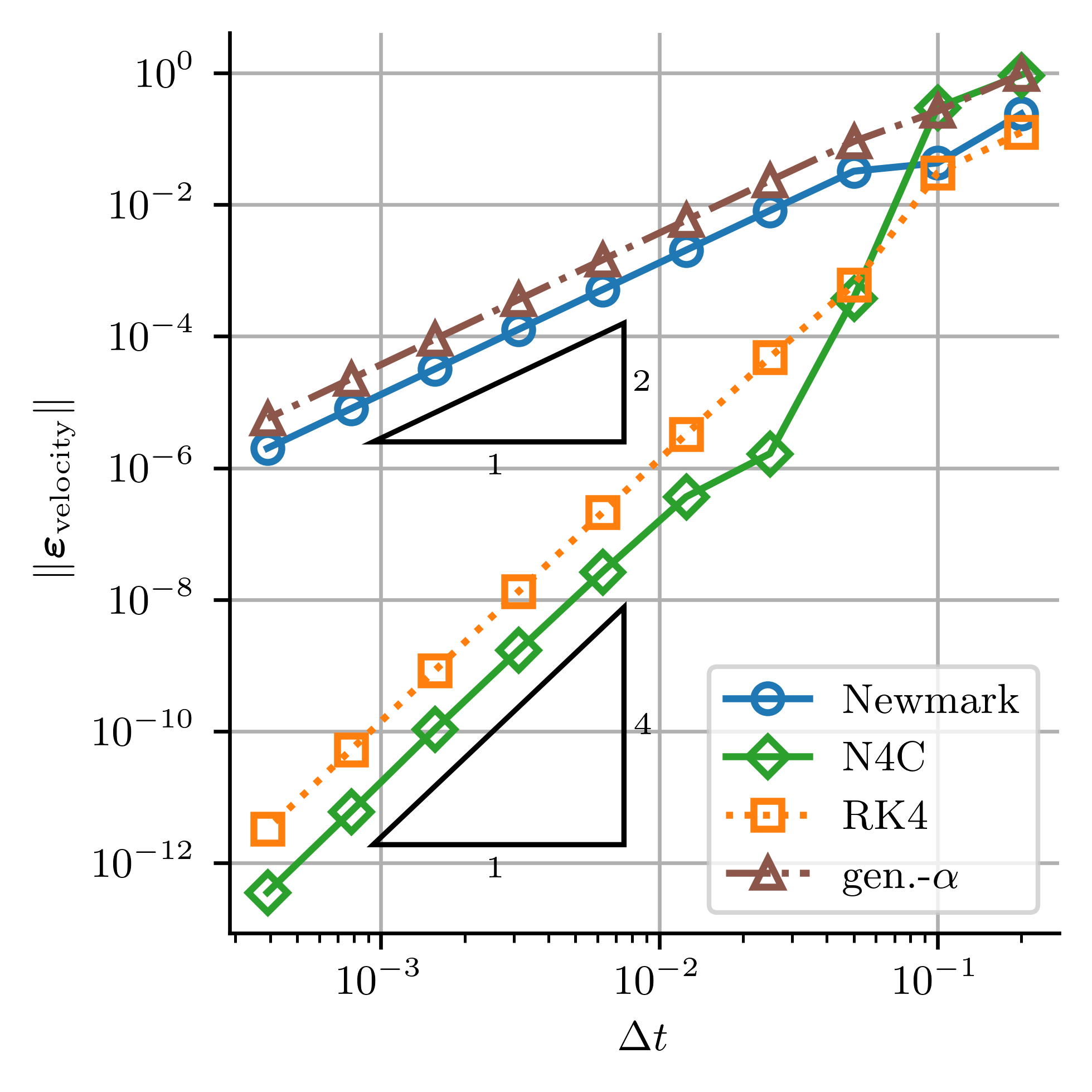}
    \caption{Convergence of the numerical solution error for the position (left) and velocity (right) using the Newmark scheme on the 1 DoF original system, the Newmark scheme with fourth-order compensation (N4C), a fourth-order Runge--Kutta scheme (RK4) and the generalized-$\alpha$ method.}
    \label{fig:comp4_error_qv}
\end{figure}

As in the case of the first convergence study (Section~\ref{sec:num_verification}), the simulations were ran using progressively finer time steps $\Dt$, and were all evaluated at $t = 0.4$. For the generalized-$\alpha$ method, here and in the following $\rho_{\infty}=1$ has been used, as this value yields the closest equivalent of the Newmark parameters used. Fig.~\ref{fig:comp4_error_qv} shows the convergence of the numerical position and velocity solutions, respectively. It is clearly visible from these graphs that the Newmark scheme with the fourth-order compensated system indeed shows a fourth-order accuracy as a function of the time step. Meanwhile, the same Newmark simulation with the exact same parameters and time step only yields second-order convergence, as does the generalized-$\alpha$ method.

% Two column figure

Furthermore, in this example, the fourth-order compensated Newmark simulation consistently gives a higher accuracy compared to an RK4 simulation with the same (sufficiently small) time step. This can be attributed to the fact that while the Runge--Kutta family of schemes are general methods, the Newmark scheme is designed especially for the simulation of second-order systems in the form of~\eqref{eq:structural} -- i.e., it can be seen as belonging to the class of structure-preserving numerical schemes. 

It is worth emphasizing here that after calculating the compensated matrices $\Chm$ and $\Khm$ according to \eqref{eq:Chm}--\eqref{eq:Khm} as well as the compensated forcing $\Fhvft$ according to \eqref{eq:Fhvft}, the fourth-order accurate calculations need the same amount of computing power as the original second-order Newmark scheme in case of direct integration. Thus, either larger time steps can be used to achieve the same accuracy, or higher accuracy can be reached using the same time step size.

After this elementary example, we turn to a higher DoF system.

\subsubsection{Numerical demonstration (3 DoF, damped, harmonic excitation)}\label{sec:comp4_3dof_damped_harmonic}
For the next demonstration of the fourth-order compensation, the system described in Section~\ref{sec:num_verification} is used, with the same time step, but with Newmark parameters $\gamma=1/2$ and $\beta=1/6$ as required for this compensation, and for the generalized-$\alpha$ method with $\rho_{\infty}=1$ accordingly.

\begin{figure}[h]
    \centering
    \includegraphics[width=\textwidth]{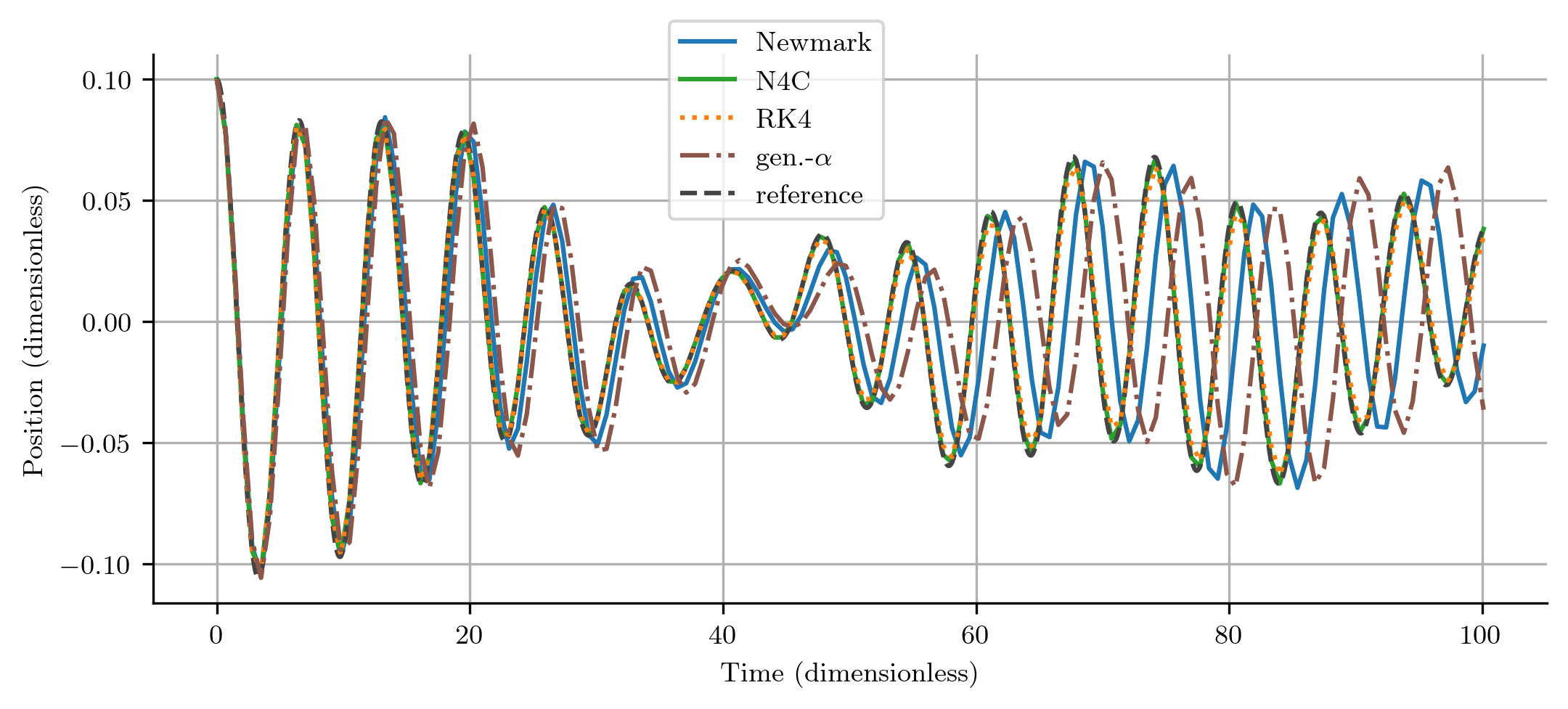}
    %\caption{Solution for the first component of the position vector $\qv$ as the function of time, from simulating the original, damped system using the Newmark method, compared to the Newmark fourth-order compensated (N4C) system also simulated with the Newmark method, against a Runge--Kutta (RK4) and a reference solution. (Time step for the former three was identically $\Dt=0.7$.)}
    \caption{Solution for the first component of the position vector as the function of time for the damped 3 DoF system with harmonic excitation. Newmark fourth-order compensation (N4C) against uncompensated Newmark, RK4, generalized-$\alpha$ and reference solutions.}
    \label{fig:comp4_cos_position_time_t100}
\end{figure}

\begin{figure}[h]
    \centering
    \includegraphics[width=0.5\textwidth]{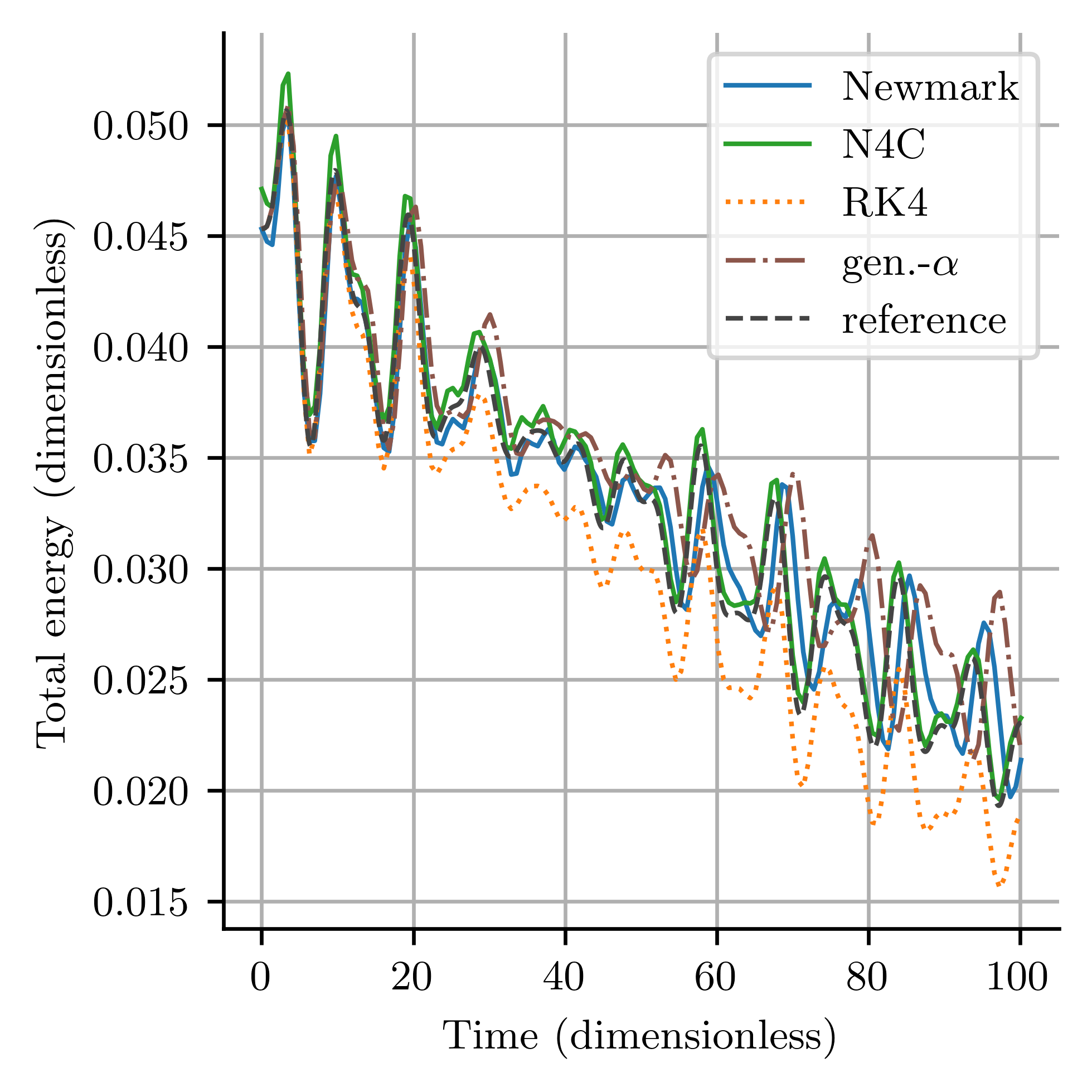}
    %\caption{Total energy as the function of time, from simulating the original, damped system using the Newmark method, compared to the Newmark damping compensated (NDC) system also simulated with the Newmark method, against a Runge--Kutta (RK4) and a reference solution. (Time step for the former three was identically $\Dt=0.7$.)}
    \caption{Total energy as the function of time for the damped 3 DoF system with harmonic excitation. Newmark fourth-order compensation (N4C) against uncompensated Newmark, RK4, generalized-$\alpha$ and reference solutions.}
    \label{fig:comp4_cos_E_tot_time_t100}
\end{figure}

Figs.~\ref{fig:comp4_cos_position_time_t100} and \ref{fig:comp4_cos_E_tot_time_t100} depict the performance of the fourth-order compensation in this scenario. In the positional solution shown in Fig.~\ref{fig:comp4_cos_position_time_t100}, the solution of the compensated system is even closer to the reference solution, compared to the results from the RK4 solver. The results deviate somewhat more if the total energy is observed in Fig.~\ref{fig:comp4_cos_E_tot_time_t100}: again, the compensated solution stays the closest to the reference solution. The uncompensated Newmark solution roughly keeps the overall damping trend but has additional inaccuracy, similarly to the generalized-$\alpha$ method. Both the Runge--Kutta and the generalized-$\alpha$ solutions also have a clear numerical dissipation error, albeit with different signs: while the RK4 method increases the observed damping, the generalized-$\alpha$ method decreases it.

% Two column figure

% Single column figure

\subsubsection{Numerical demonstration (3 DoF, undamped, non-harmonic excitation)}
\begin{figure}[h]
    \centering
    \includegraphics{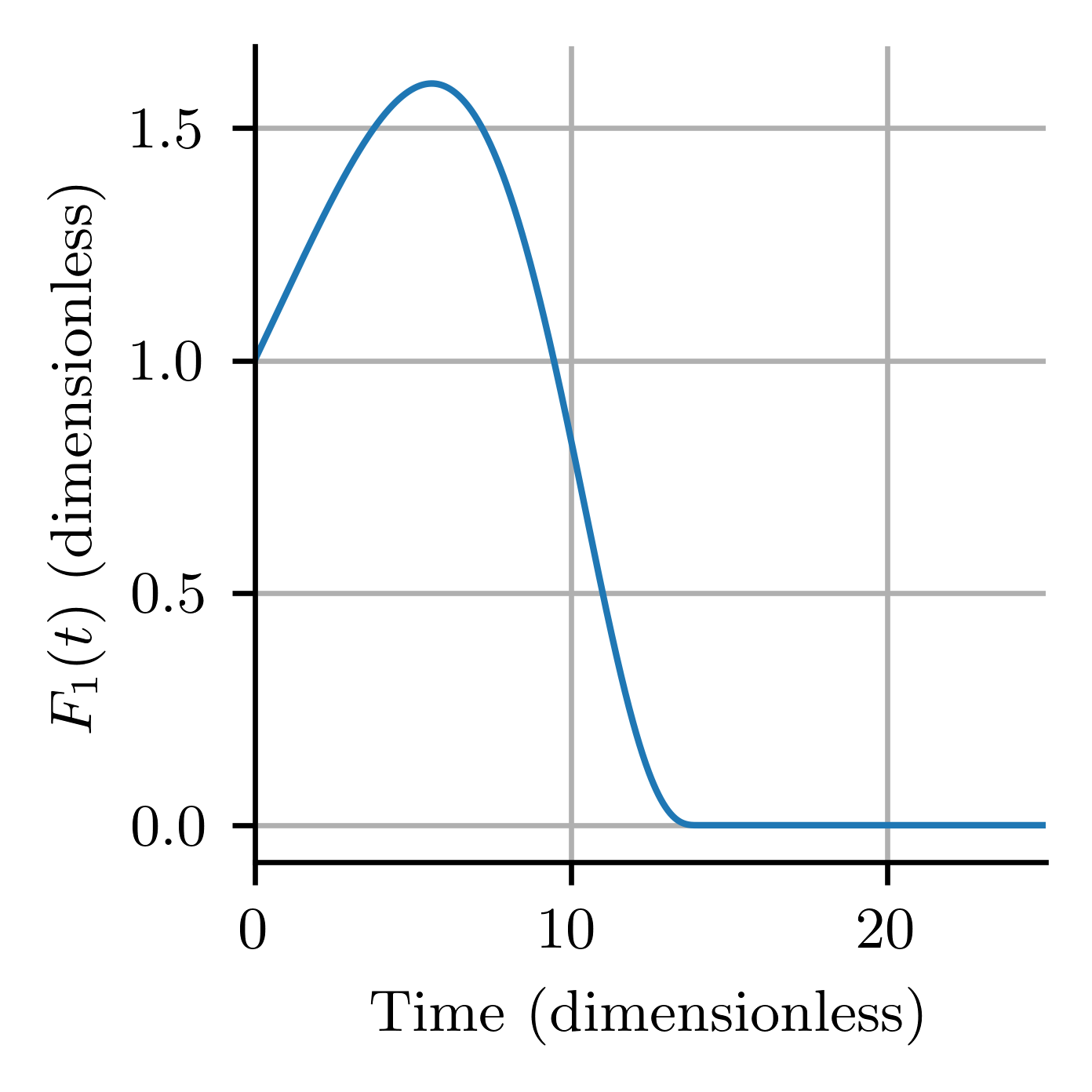}
    \caption{The first component of the excitation at the beginning of the simulation with $\mu=0.2$ and $\tstar=14$.}
    \label{fig:bump}
\end{figure}

For another demonstration of the fourth-order compensation, we also show the performance of the method under non-harmonic excitation. The excitation is a pulse, defined as
\begin{gather}
    \Fvft = 
    \begin{pmatrix}
        1
        \\
        0
        \\
        0
        \\
    \end{pmatrix}
    %\mathrm{e}^
    \exp(\frac{t}{\mu \tstar}) \left(1- \frac{t}{\tstar}\right)^{3}
    ,\,\, 0 \leq t \leq \tstar;
    \quad
    \Fvft \equiv \mx{0},\, t > \tstar,
\end{gather}
where $\mu$ is a shape parameter and $\tstar$ is the cutoff time for the pulse. The shape of the pulse is shown in Fig.~\ref{fig:bump} for the values of these parameters used in the simulations.

% Single column figure

\begin{figure}[h]
    \centering
    \includegraphics[width=\textwidth]{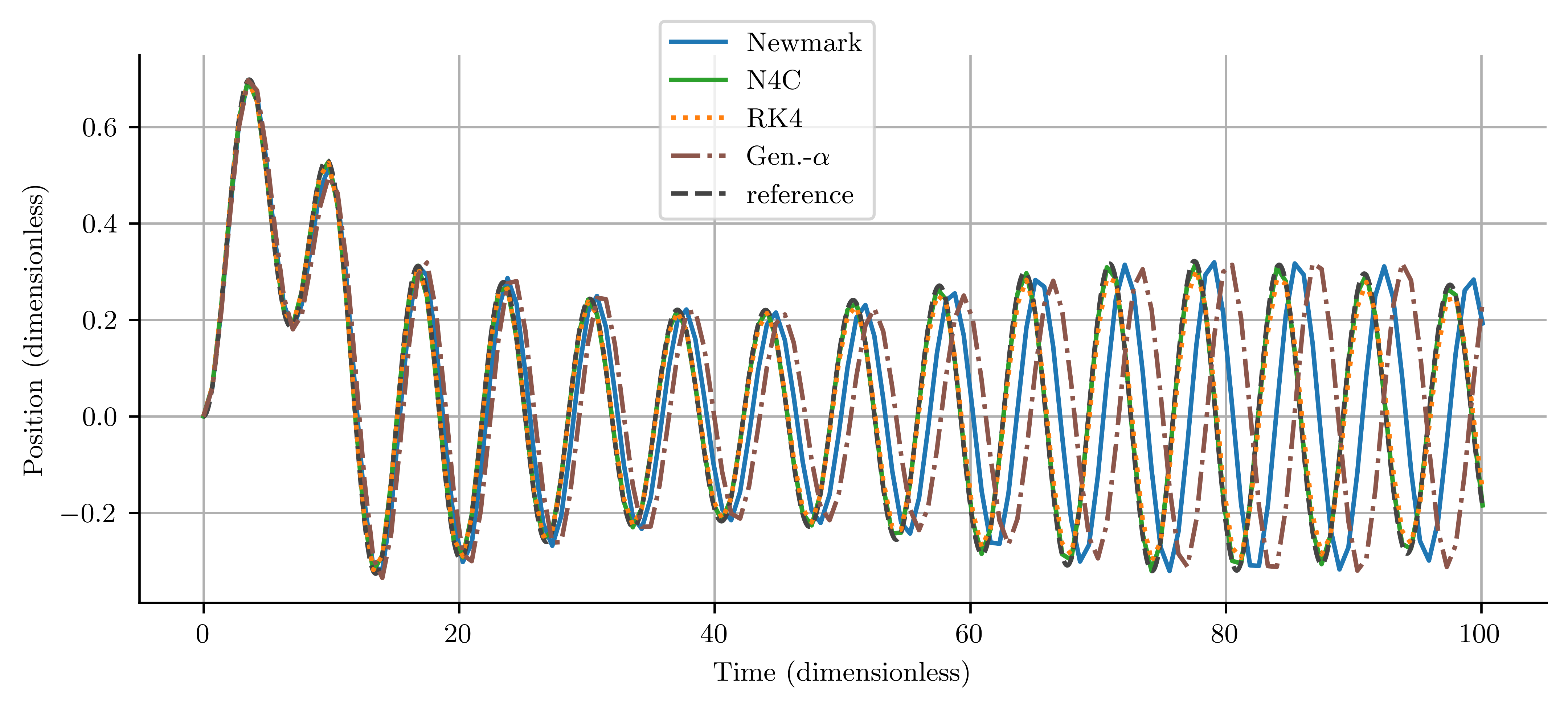}
    \caption{Solution for the first component of the position vector as the function of time for the undamped 3 DoF system with finite pulse excitation. Newmark fourth-order compensation (N4C) against uncompensated Newmark, RK4, generalized-$\alpha$ and reference solutions.}
    \label{fig:comp4_bump_position_time_t100}
\end{figure}

\begin{figure}[h]
    \centering
    \includegraphics[width=0.5\textwidth]{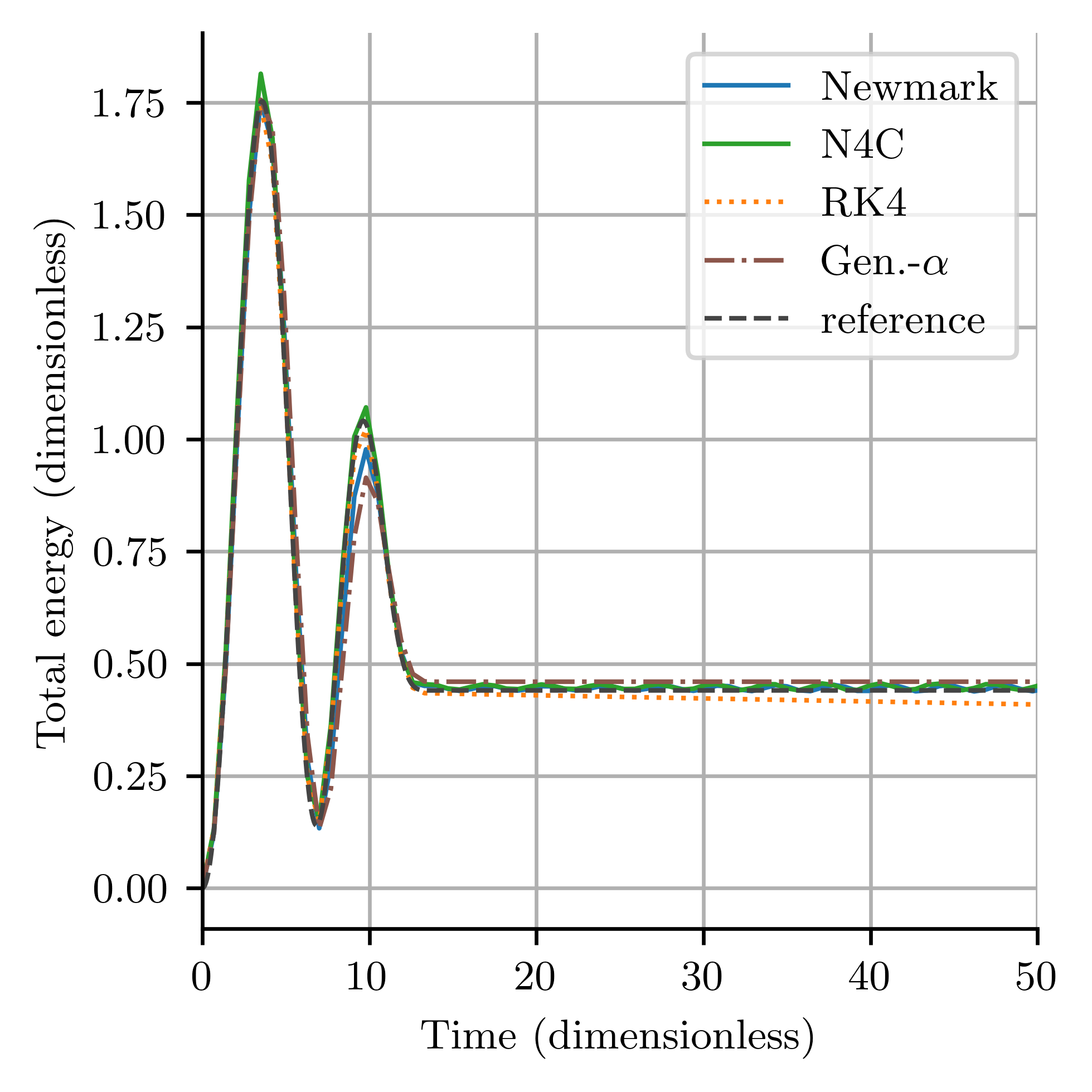}
    \caption{Total energy as the function of time for the damped 3 DoF system with finite pulse excitation. Newmark fourth-order compensation (N4C) against uncompensated Newmark, RK4, generalized-$\alpha$ and reference solutions.}
    \label{fig:comp4_bump_E_tot_time_t100}
\end{figure}

As the excitation has been chosen to introduce a finite amount of energy, setting $\Chm=\mx{0}$ allows the energy preservation of the fourth-order compensation to be testable for $t>\tstar$. Figs.~\ref{fig:comp4_bump_position_time_t100} and~\ref{fig:comp4_bump_E_tot_time_t100} show the position and total energy solution as the function of time, respectively. Once again, the compensated solution reproduces the position the most faithfully over time, and its energy preserving behaviour is superior to that of the Runge--Kutta method, while being similar to an uncompensated Newmark or generalized-$\alpha$ simulation with the same settings.

% Two column figure

% Single column figure

\subsubsection{Numerical demonstration (non-differentiable excitation, numerical derivatives)}

% During the derivations of the fourth-order compensation, due to the mathematical requirements of BEA, it has been assumed that the excitation is continuous and is differentiable at least twice, as \eqref{eq:Fhvft} contains the second derivative of $\Fvft$. However, in practical applications, the excitation available does not always fulfil these requirements: it is possible that the excitation is only available from measurements at discrete time instants, and these instants can even differ from the time steps used in the simulation. At first glance, this seems to significantly limit the applicability of our method for non-continuous excitations.
% 
% Nevertheless, our experience is that the compensation \eqref{eq:Fhvft} -- derived for continuous functions -- can still be approximated using appropriately chosen numerical formulas for derivatives, making the fourth-order compensation applicable to non-continuous, non-differentiable excitations as well.

During the derivations of the fourth-order compensation, due to the mathematical requirements of BEA, it has been assumed that the excitation is continuous and is differentiable at least twice, as \eqref{eq:Fhvft} contains the second derivative of $\Fvft$. In practical applications, the excitation available does not always fulfil these requirements: it is possible that the excitation is only available from measurements at discrete time instants, and these instants can even differ from the time steps used in the simulation; or the excitation might not even be differentiable. At first glance, these may seem to limit the applicability of our method for non-continuous or non-differentiable excitations.

Nevertheless, our experience is that the compensation \eqref{eq:Fhvft} -- derived for continuous functions -- can still be approximated using appropriately chosen numerical formulas for derivatives, making the fourth-order compensation applicable to non-continuous, non-differentiable excitations as well.

\begin{figure}[hbt]
    \centering
    \includegraphics[width=0.495\textwidth]{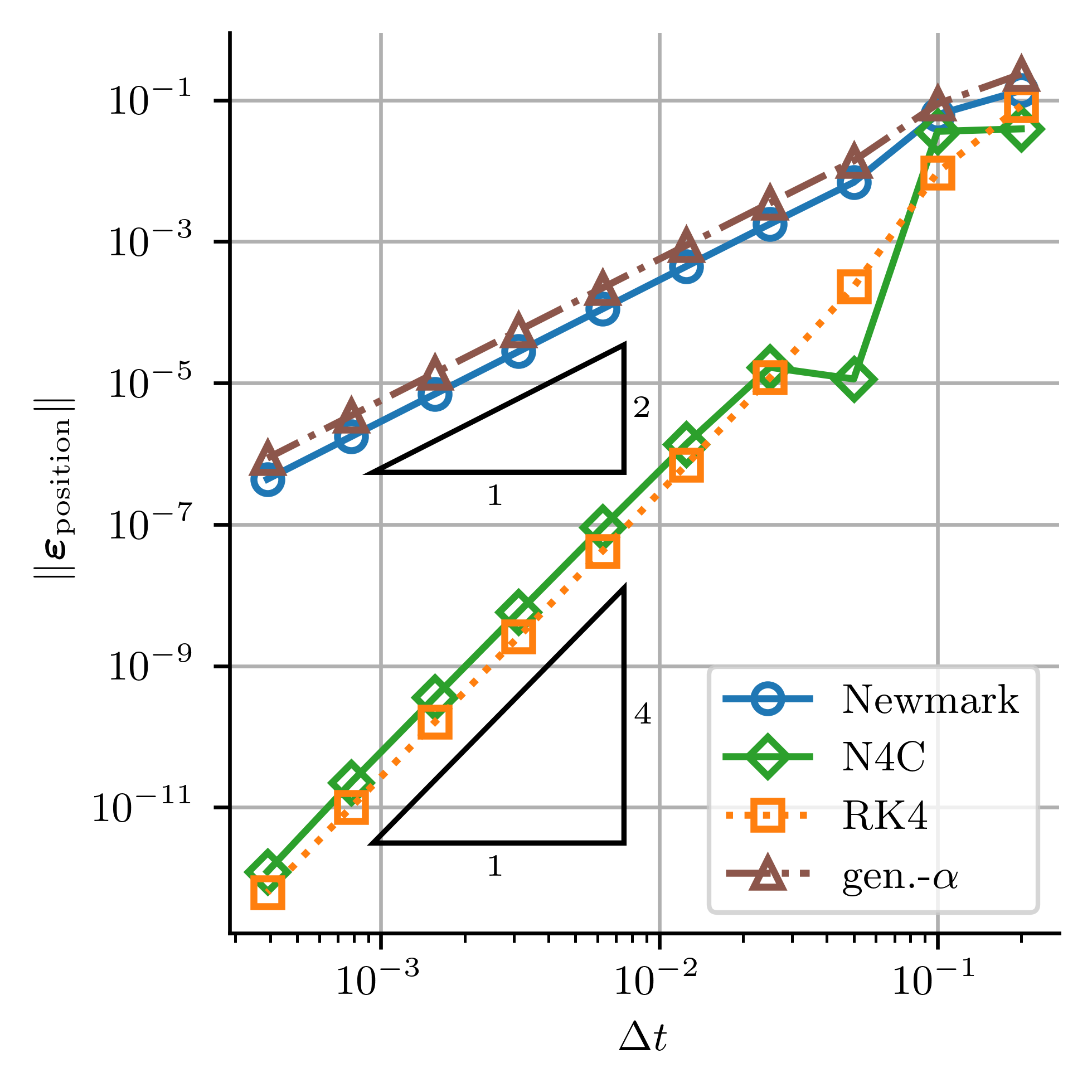}
    \includegraphics[width=0.495\textwidth]{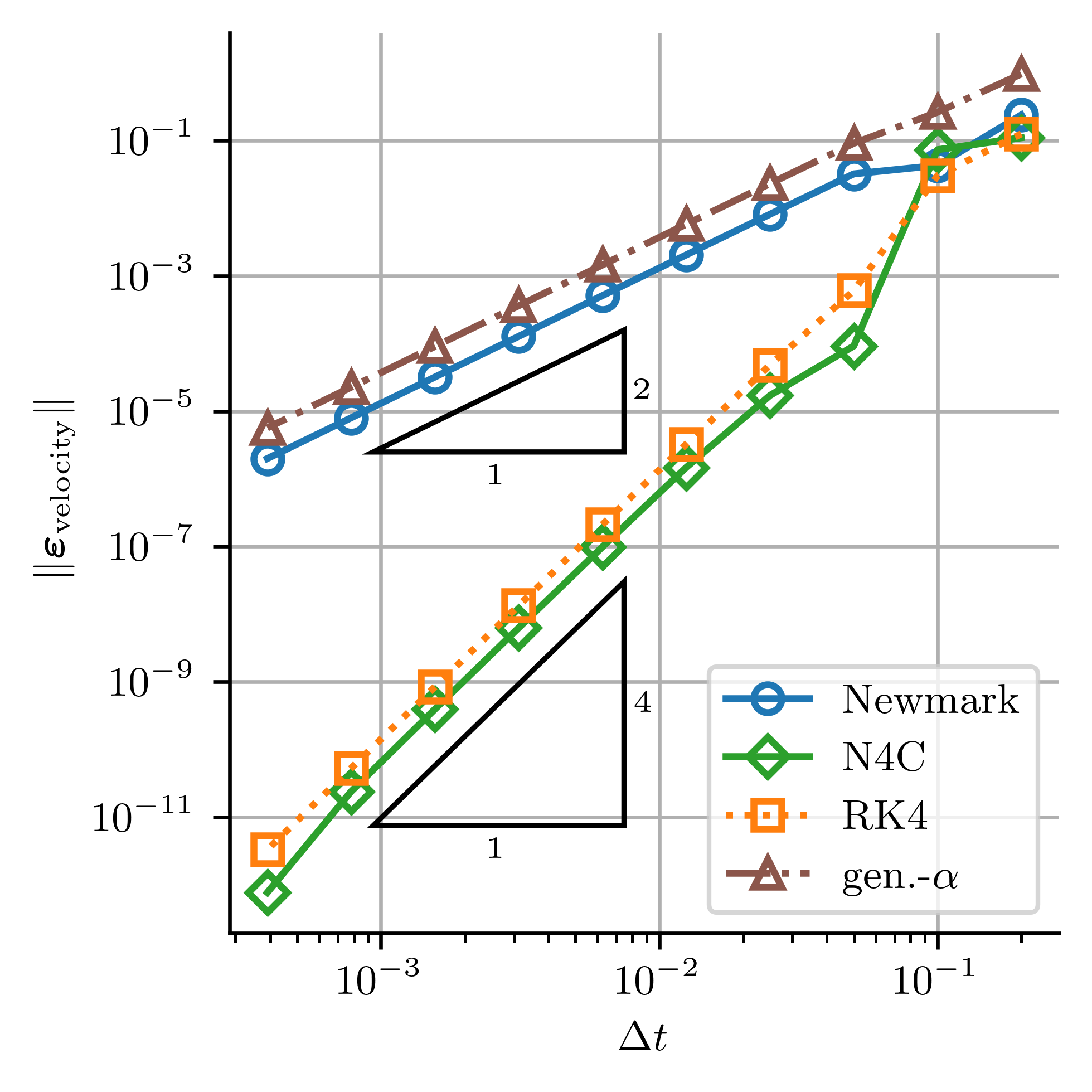}
    \caption{Convergence of the numerical solution error for the position (left) and velocity (right) using the Newmark scheme on the 1 DoF original system, the Newmark scheme with fourth-order compensation with numerically approximated derivatives (N4C), a fourth-order Runge--Kutta scheme (RK4) and the generalized-$\alpha$ method.}
    \label{fig:comp4_num_error_qv}
\end{figure}

Equation \eqref{eq:Fhvft} shows that the derivatives $\Fvftp$ and $\Fvftpp$ are in the $\Ord\left({\Dt^2}\right)$ terms, thus the use of second-order accurate numerical differentiation formulas are appropriate: the error introduced by them will be $\Ord\left({\Dt^4}\right)$ with respect to the entire equation. Subsequently, we can use
\begin{align}
    % Fv_t_num_ = (F_(t_ + h_) - F_(t_ - h_))/(2*h_)
    % Fvv_t_num_ = (F_(t_ + h_) - 2*F_(t_) + F_(t_ - h_))/(h_**2)
    \Fvftp &\approx \frac{\Fv\leftf(t + \Dt \rightf) - \Fv\leftf(t - \Dt \rightf)}{2 \Dt}, \label{eq:Fvftpapprox}
    \\
    \Fvftpp &\approx \frac{\Fv\leftf(t + \Dt \rightf)- 2 \Fv\leftf(t \rightf) + \Fv\leftf(t - \Dt \rightf)}{\Dt^2}, \label{eq:Fvftppapprox}
\end{align}
where $\Fvft$ is either available (without its derivatives) at $t$ or is formulated using a suitably high-degree interpolation of measurement data. Using the above \eqref{eq:Fvftpapprox}--\eqref{eq:Fvftppapprox}, we can repeat the convergence analysis detailed in Section~\ref{sec:convergence1dof}. The results are shown in Fig.~\ref{fig:comp4_num_error_qv}, demonstrating that though the overall error increases slightly, the numerically approximated derivatives still yield a compensation that is fourth-order accurate, as expected.

% Two column figure

\begin{figure}[h]
    \centering
    \includegraphics[width=\textwidth]{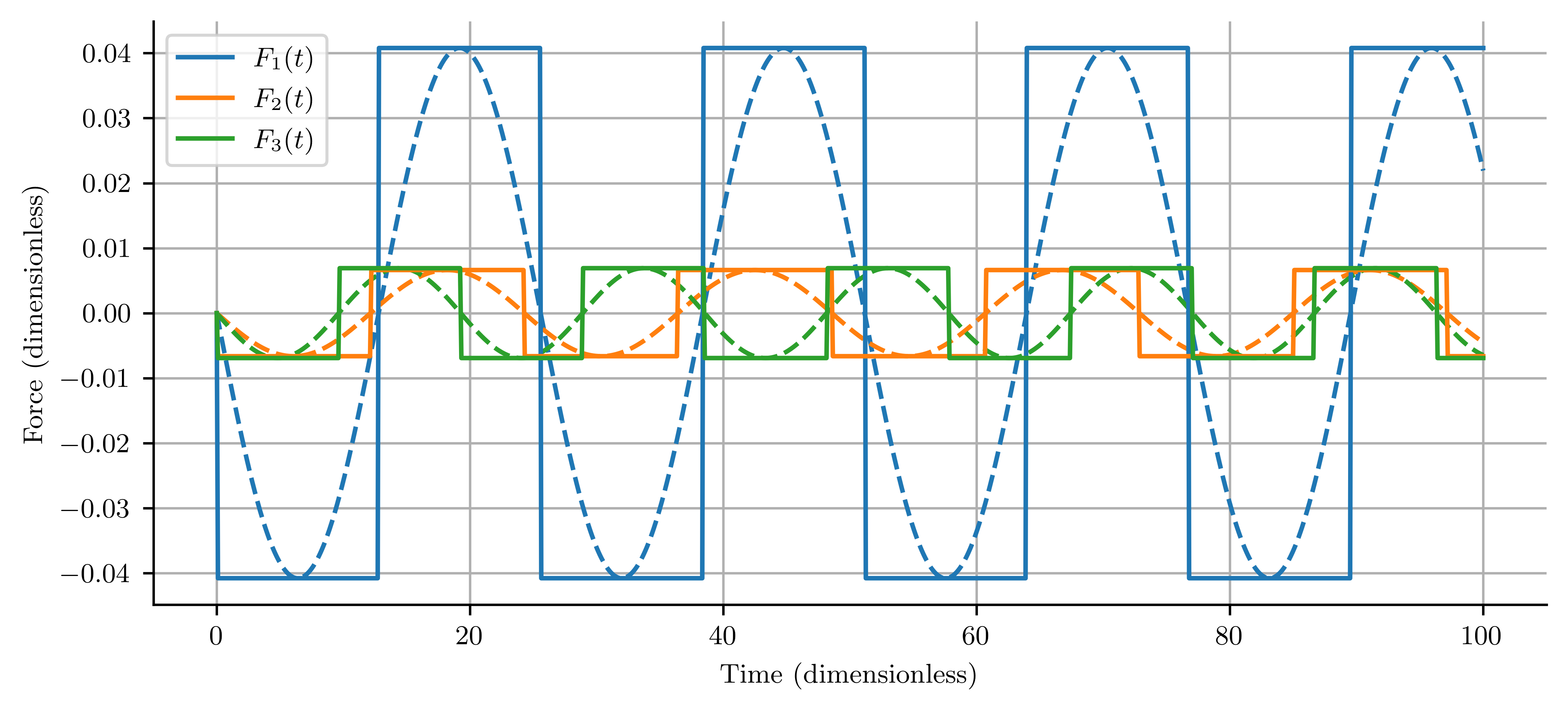}
    \caption{Square-wave excitation according to \eqref{eq:squarewave}, with the underlying sinusoidal oscillations also shown.}
    \label{fig:squarewave}
\end{figure}

\begin{figure}[h]
    \centering
    \includegraphics[width=\textwidth]{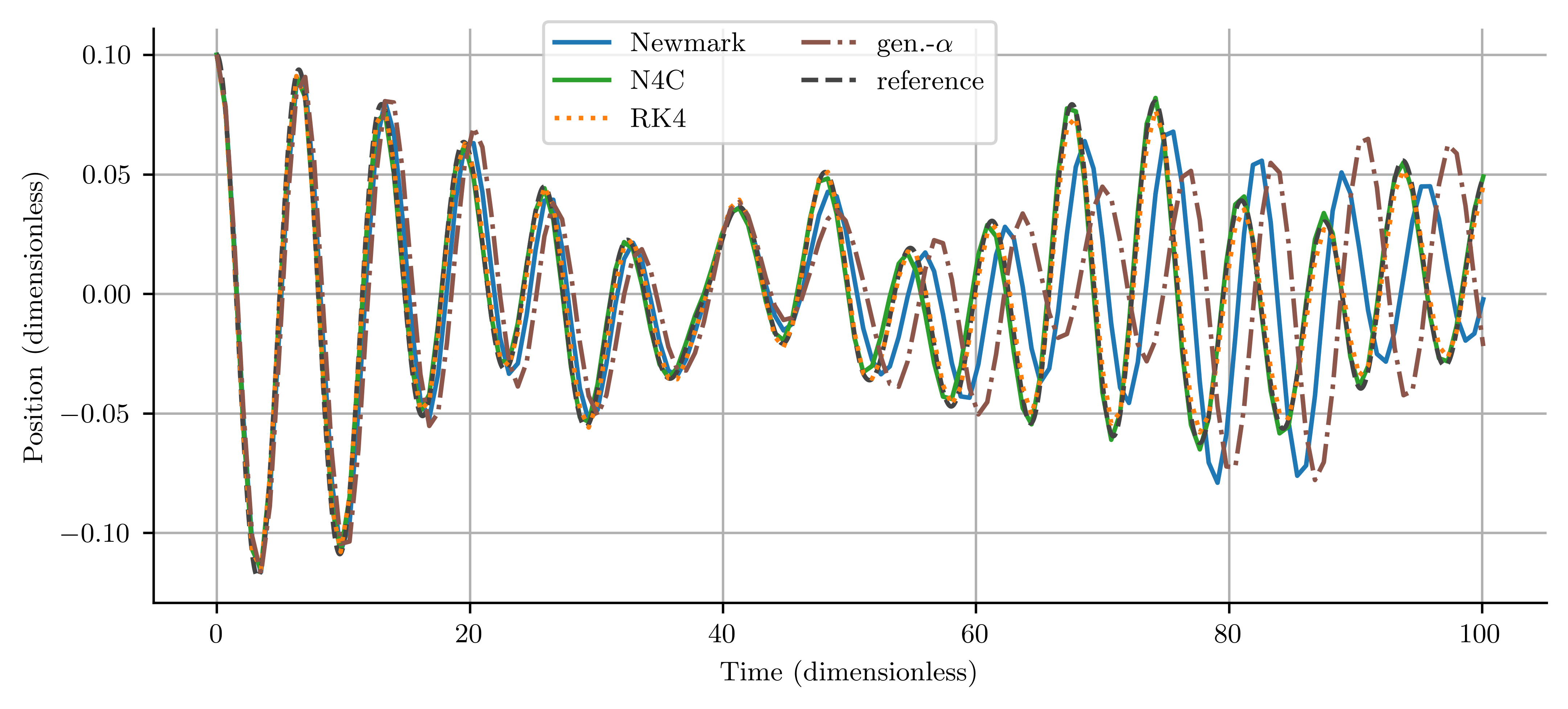}
    \caption{Solution for the first component of the position vector as the function of time for the undamped 3 DoF system with square-wave excitation. Newmark fourth-order compensation with numerical derivatives (N4C) against uncompensated Newmark, RK4, generalized-$\alpha$ and reference solutions.}
    \label{fig:comp4_numdiff_square_position_time_t100}
\end{figure}

\begin{figure}[h]
    \centering
    \includegraphics[width=0.5\textwidth]{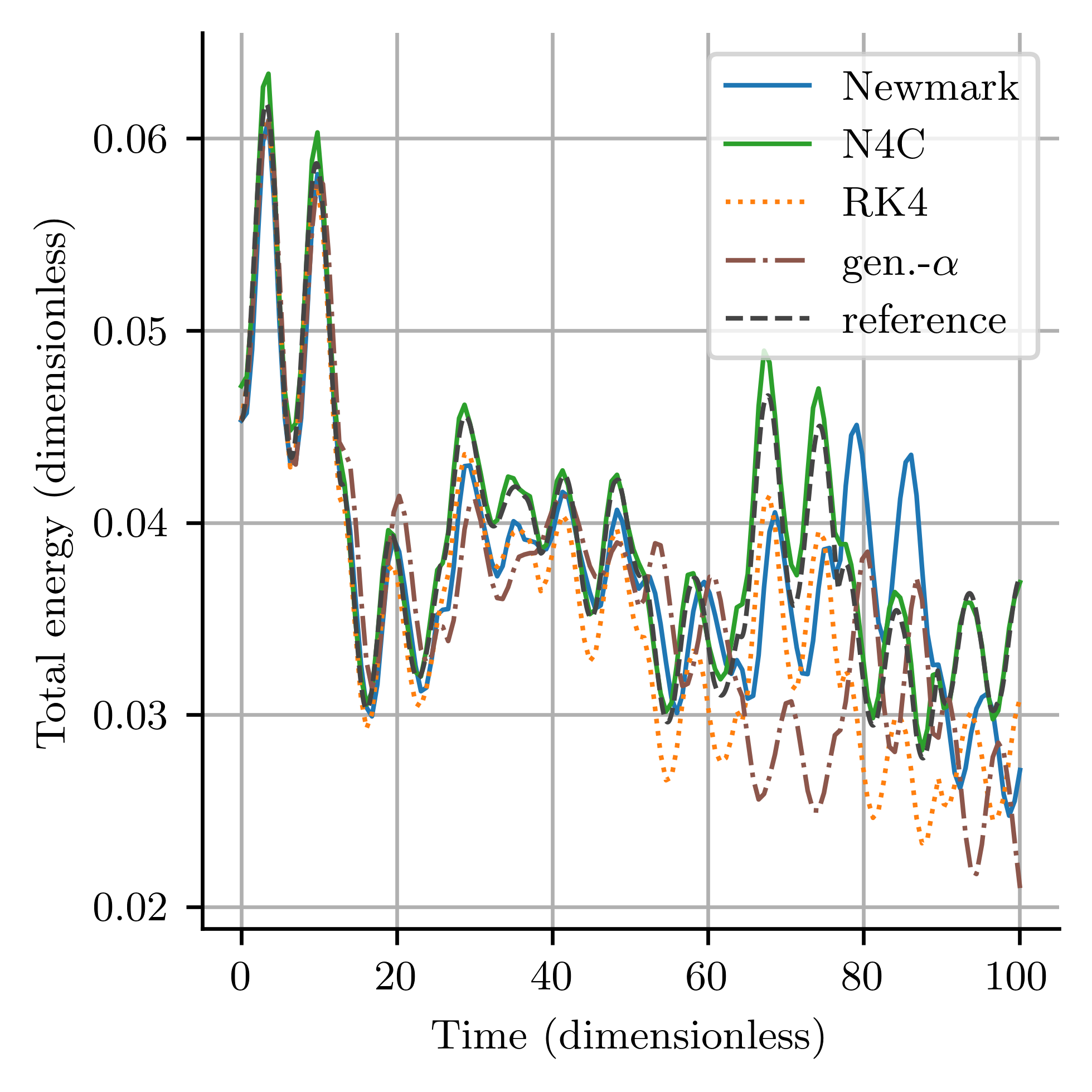}
    \caption{Total energy as the function of time for the damped 3 DoF system with square-wave excitation. Newmark fourth-order compensation with numerical derivatives (N4C) against uncompensated Newmark, RK4, generalized-$\alpha$ and reference solutions.}
    \label{fig:comp4_numdiff_square_E_tot_time_t100}
\end{figure}

As a more detailed example of using numerical approximations for non-continuous, non-differentiable excitation, we changed the sinusoidal excitation of the 3 DoF system used in Sections~\ref{sec:num_verification} and \ref{sec:comp4_3dof_damped_harmonic} to a square-wave excitation as
\begin{gather}
    \Fvft = 
    \begin{pmatrix}
        -0.040790 \cdot \sign( \sin(0.2457\cdot t) )\\ 
        -0.006630 \cdot \sign( \sin(0.2587\cdot t) )\\ 
        -0.006914 \cdot \sign( \sin(0.3262\cdot t) )\\
    \end{pmatrix},\label{eq:squarewave}
\end{gather}
also illustrated in Fig.~\ref{fig:squarewave}. The results of simulating this system are shown in Fig.~\ref{fig:comp4_numdiff_square_position_time_t100} and Fig.~\ref{fig:comp4_numdiff_square_E_tot_time_t100}. These show that not only is the fourth-order compensation with numerical derivatives applicable to non-differentiable excitations, but it also significantly outperforms the Newmark, generalized-$\alpha$ and RK4 schemes in this example. It should be also noted that even with the numerical approximation of the derivatives, the fourth-order compensated Newmark method only needs three evaluations of $\Fvft$, while the RK4 method uses four.

% Two column figure

% Two column figure

% Single column figure

\subsubsection{Numerical demonstration (high-DoF, finite element model matrices)}
\begin{figure}[h]
    \centering
    \hfill
    \includegraphics[width=0.4\textwidth]{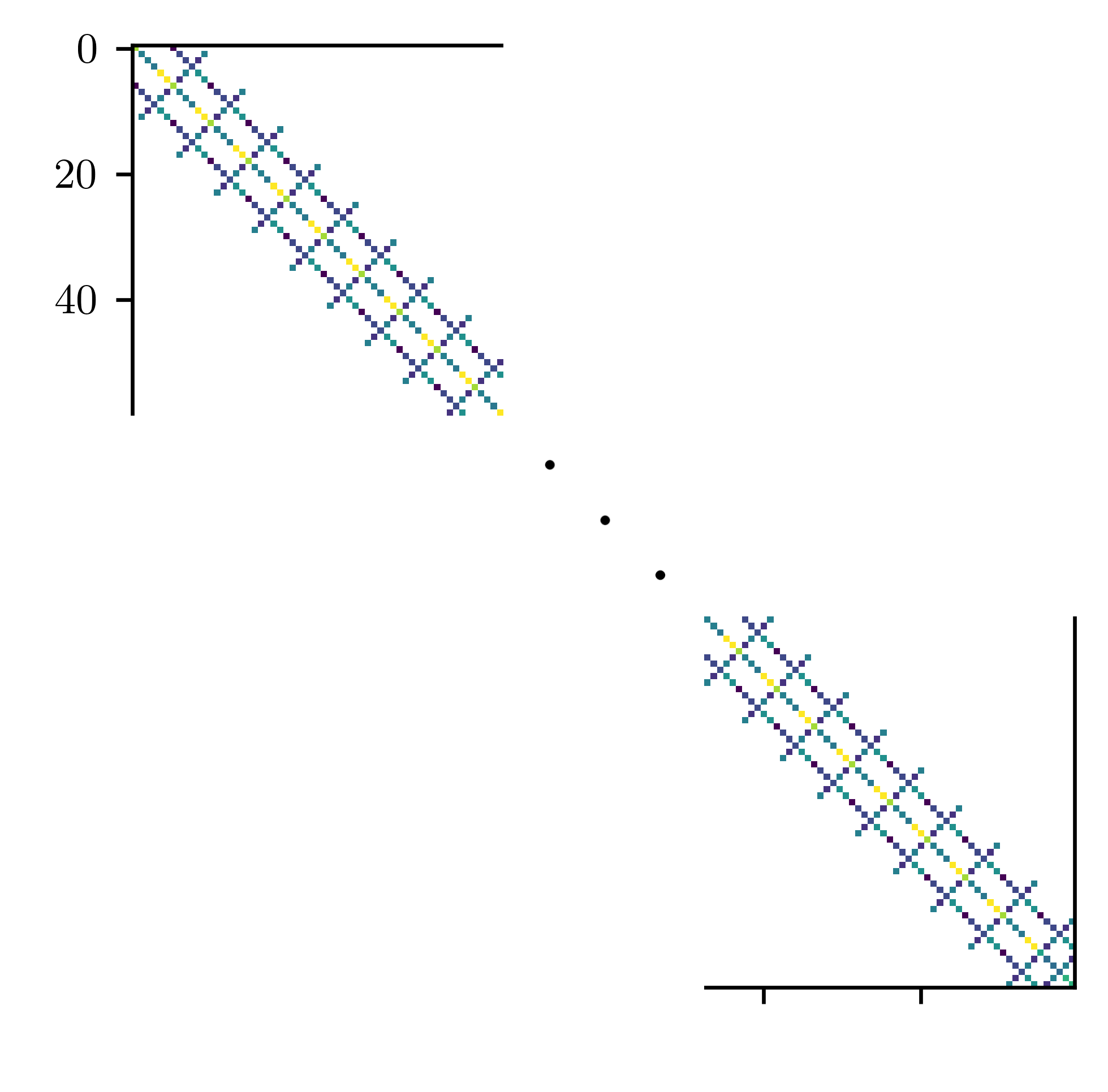}
    \hfill
    \includegraphics[width=0.4\textwidth]{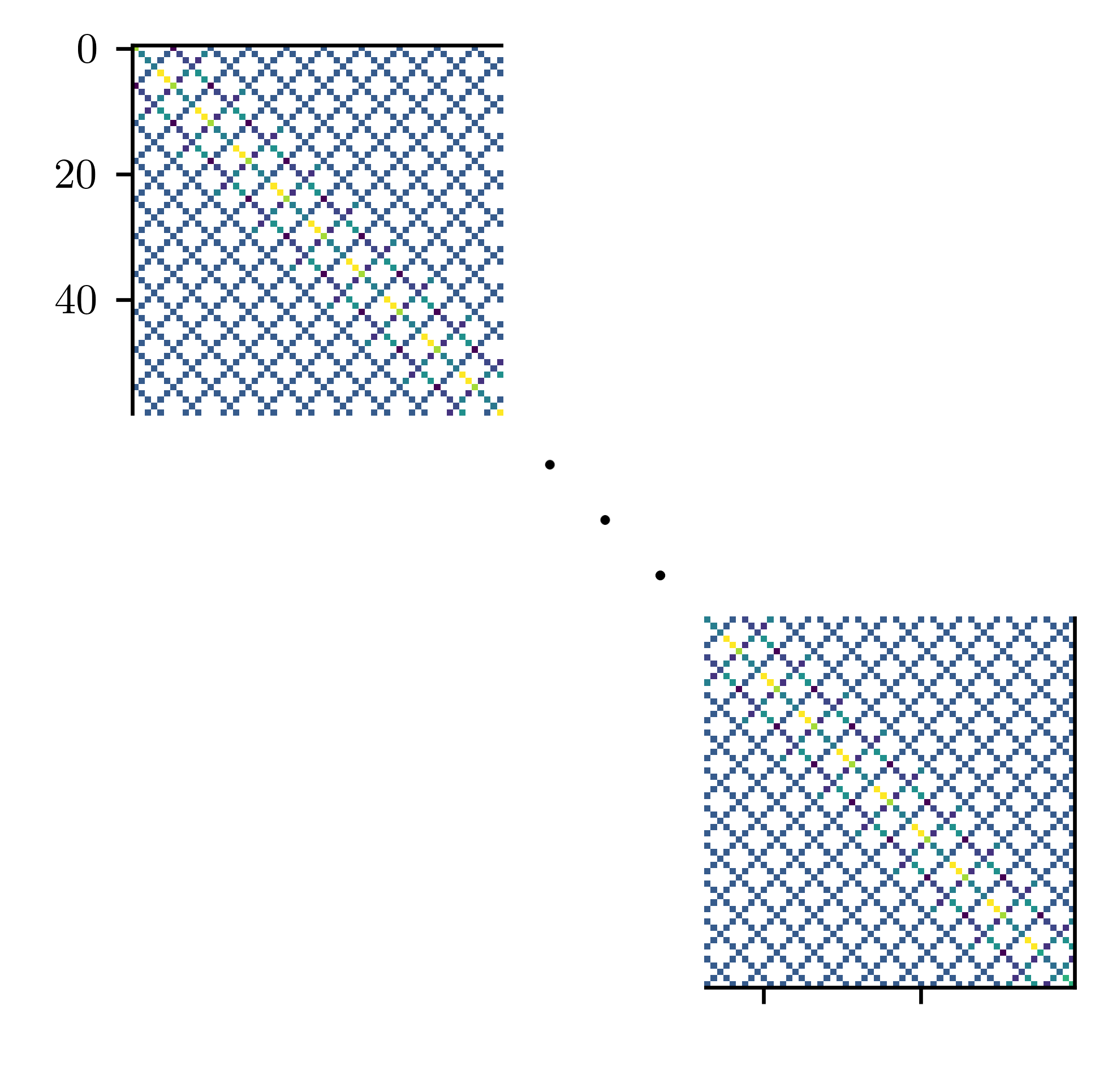}
    \hfill
    \caption{Sparsity structure of the original finite element stiffness matrix $\Km$ (left) versus the compensated stiffness matrix $\Khm$ (right).}
    \label{fig:sparsity}
\end{figure}

So far, the examples shown have all been of systems with at most 3 DoF, and having dense matrices. However, the Newmark method and its variants are also widely used for solving finite element problems which have a high DoF and matrices with a sparse structure which is also exploited for efficient calculations. On the other hand, the compensated matrices \eqref{eq:Chm}--\eqref{eq:Khm} will become dense, as illustrated in Fig.~\ref{fig:sparsity}. This raises the question whether the computational disadvantages posed by non-sparse matrices outweigh the advantages of fourth-order convergence. In our last example, we explore the characteristics of the fourth-order compensation in this context.

% Two column figure

\begin{figure}[h]
    \centering
    \includegraphics[width=0.5\textwidth]{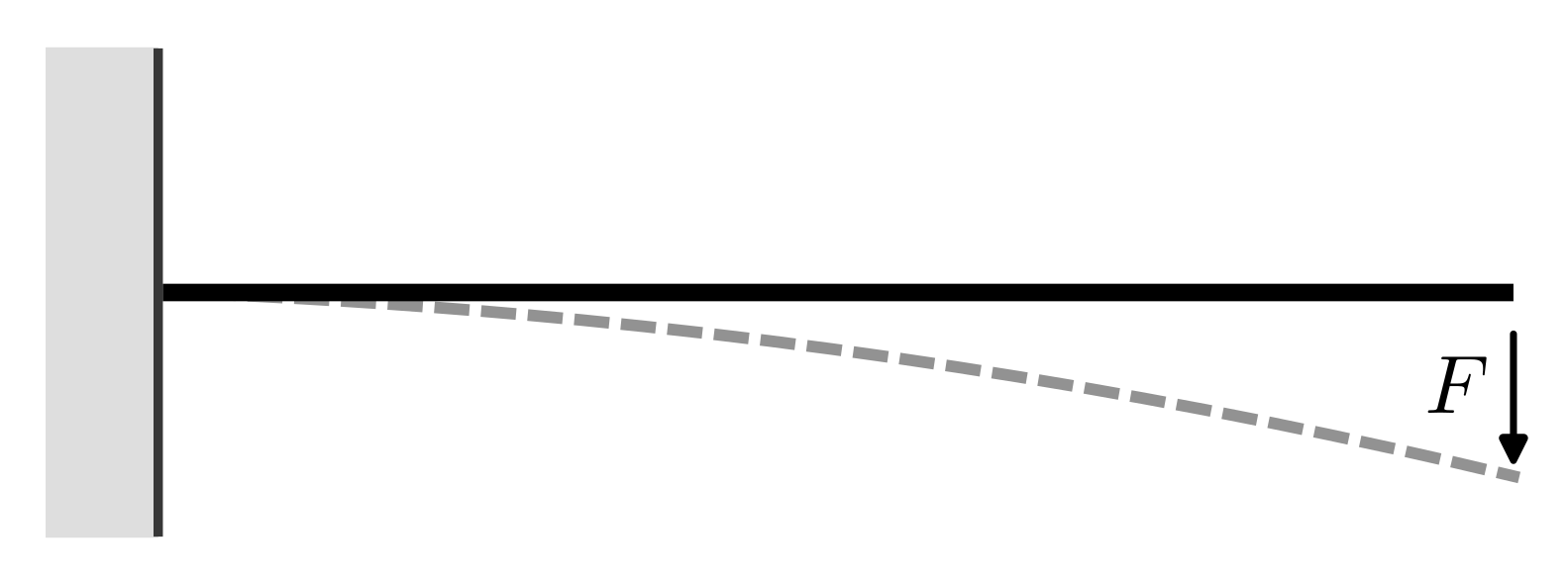}
    \caption{Initial static loading of cantilever beam modelled. (Not to scale.)}
    \label{fig:cantilever_loading}
\end{figure}

We consider a one-meter steel cantilever beam with a diameter of 10 mm, modelled in ANSYS Mechanical using 250 elements of type BEAM188 \cite{ansys2023elements}. As BEAM188 elements have two nodes with 6 DoF each, this results in a system with 756 DoF, which has been used directly after being exported from ANSYS. The matrix structure (shown in Fig.~\ref{fig:sparsity}) is very sparse, even with finite element standards, due to the geometry. The simulation is a free vibration with an initial condition from a static deformed state due to a vertical loading of 10 N at one end (as shown in Fig.~\ref{fig:cantilever_loading}).

% Single column figure

For comparison, the fourth-order compensated Newmark method is tested against the generalized-$\alpha$ method. (The non-compensated Newmark method has been consistently found to be of lower accuracy than the generalized-$\alpha$ method, with no significant advantages in computing time, thus it has been excluded from this comparison for clarity.) To contrast the performance  of the two methods, the same sparse-matrix solver is used: as the Newmark method is a special case of the generalized-$\alpha$ method, using $\alpha_f=\alpha_m=0$ for the former makes a direct computational comparison with the latter possible. For the generalized alpha method, $\rho_{\infty}=1$ is used.

For the different simulations, the time step has been halved from the maximum value allowable by the stability criterion until the convergence of the error is inhibited by the limits of numerical precision. For each time step size, the simulations were ran five times to get a more accurate estimate of the runtime; the maximum runtime has been used to normalize the computing time shown. (Similarly to previous examples, an RK4 solution with significantly finer time steps has been used as a reference.)

\begin{figure}[h]
    \centering
    \includegraphics[width=0.495\textwidth]{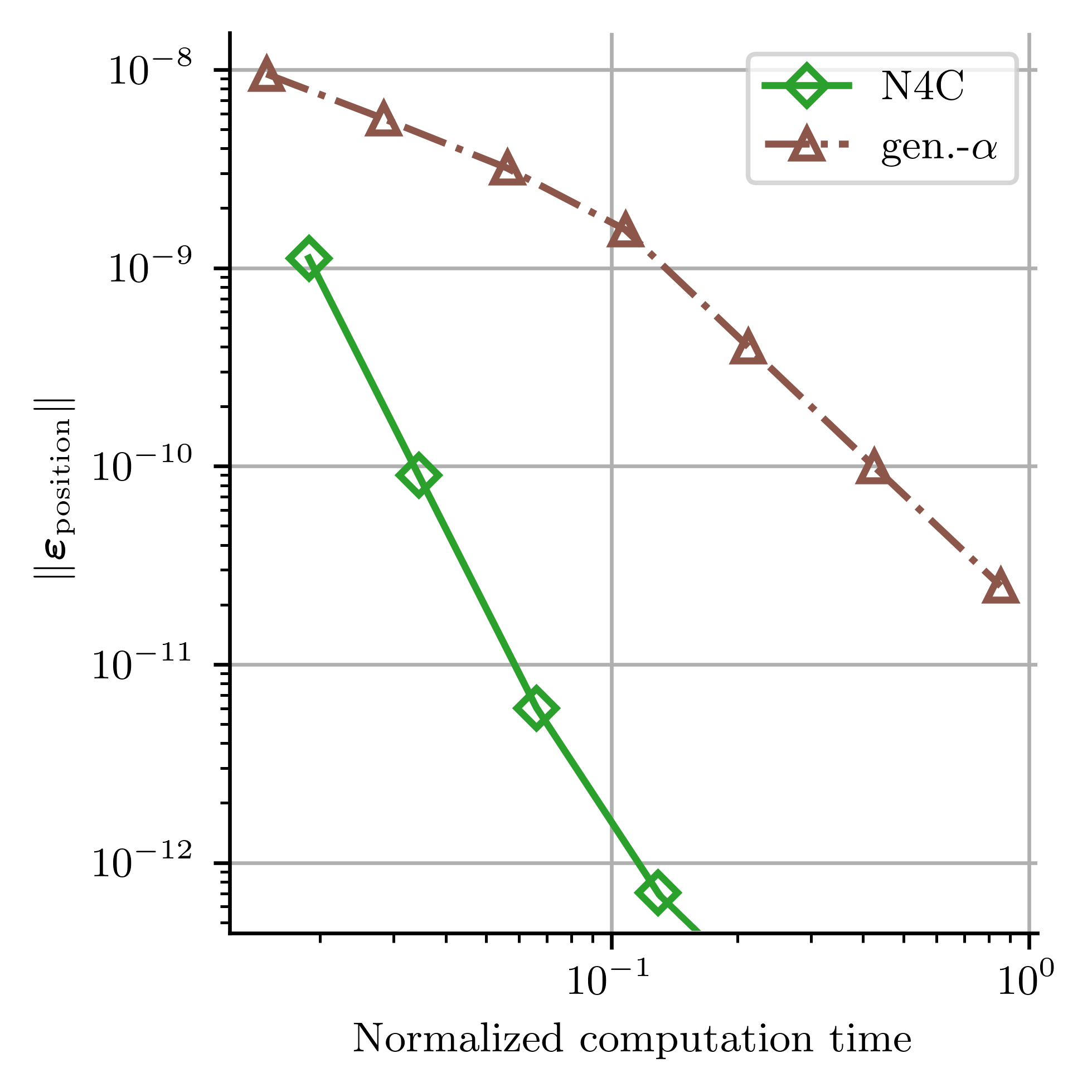}
    \includegraphics[width=0.495\textwidth]{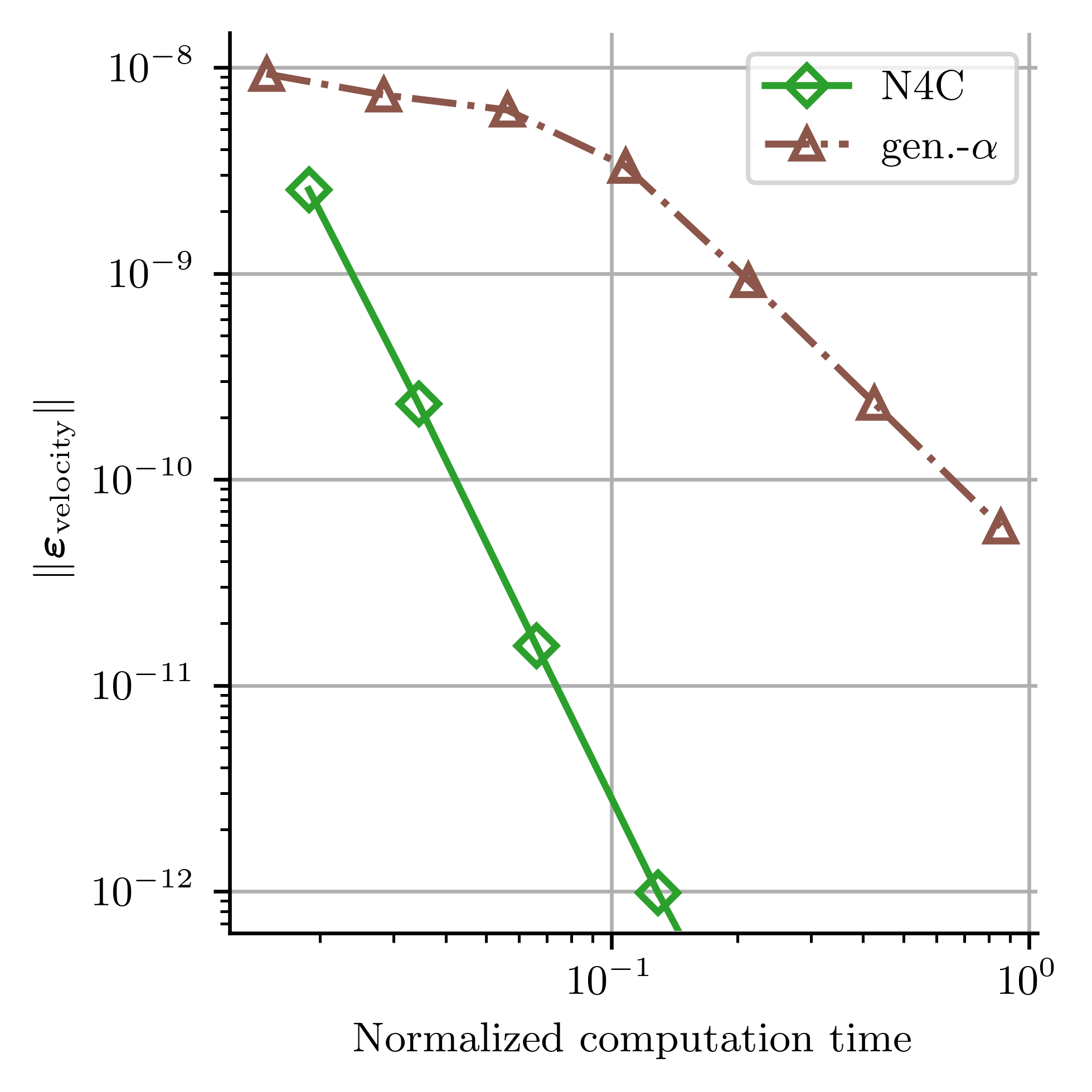}
    \caption{Accuracy of simulations against computing time, for the high-DoF finite element beam problem. Fourth-order compensated Newmark (N4C) against the generalized-$\alpha$ method.}
    \label{fig:fea_convergence}
\end{figure}

Fig.~\ref{fig:fea_convergence} shows the accuracy of the simulations versus the computing time needed to achieve the results. It is clear that even though the sparse matrices of the generalized-$\alpha$ method yield faster computations for the same time step size, the advantages of fourth-order convergence result in significantly better accuracy for the compensated Newmark method. Equivalently, the same accuracy can be achieved using the fourth-order compensated Newmark method with significantly larger timesteps, and thus significantly lower computation time. We have found that this trend persists for even larger, several thousand DoF discretizations of the same problem. A more in-depth comparison for even larger systems with more complex geometry, and a detailed analysis of memory use of the two algorithms is left for future work.

% Two column figure

It is to be stressed at this point that practical implementations of the Newmark and other methods usually circumvent the inversion of the mass matrix $\Mm$, especially for non-lumped-mass element formulations. This approach is also included in the generalized-$\alpha$ solver used for the above calculations, and the calculation of our compensated formulas \eqref{eq:Chm}--\eqref{eq:Fhvft} also allows for such an implementation at no additional cost. More specifically, $\Mminv$ only occurs in combinations such as $\Mminv \Km$, $\Mminv \Cm$, $\Mminv \Fvft$ etc., and during the above comparisons, we calculated these directly instead of multiplying by $\Mminv$.

As a future possibility for enhancing our method, the compensated matrices may be made even more computationally efficient if we relax the self-imposed constraint of \eqref{eq:newmark_mod_secondorder}, where we have fixed an ambiguity by keeping the original mass matrix in the distorted and compensated systems. This equation being definite only up to (maximal-rank) matrix multiples allows for further optimizations in practical applications. Finding such a computationally more efficient variant of \eqref{eq:newmark_mod_secondorder} would be worthy of future investigation.

\section{Discussion}\label{sec:discussion}
Through the use of backward error analysis, we have derived the modified or distorted equation corresponding to the Newmark method applied to linear structural systems under transient excitation. We have also formulated the distorted equation in terms of a second-order equation comparable to the original system containing distorted damping, stiffness and excitation terms. Both of these results are valid up to the order of $\Dt^2$, while additional terms can be calculated using the supplied computer algebra code, if necessary.

We have used the distorted equations to derive two compensation strategies to improve the numerical results calculated using the Newmark scheme. The first eliminates the numerical damping introduced by the Newmark scheme through a compensated damping matrix for both damped and undamped linear systems. The compensation is valid for arbitrary values of the Newmark parameters $\gamma$ and $\beta$. Naturally, there exist cases where the numerical damping effect is desirable \cite{serfozo2023method}: thus, in a similar vein, any desired numerical damping profile that can be expressed explicitly as an additional part of the original damping matrix could be achieved. Further research
%is needed
in this direction
%to
may more closely
specify the type of numerical damping that could be achieved through this method.

The second compensation strategy achieves fourth-order accuracy in the originally second-order Newmark method. This approach employs compensated damping and stiffness matrices, as well as a compensated excitation term, and is valid for the specific parameter combination of $\gamma=1/2,\,\beta=1/6$.
%Further investigations
%are needed to determine whether
%may be able to explore how well
%the improved accuracy persists for non-smooth excitation functions, as the compensated excitation $\Fhvft$ contains derivatives of $\Fvft$. 
Even though the compensated excitation $\Fhvft$ contains derivatives of $\Fvft$, we have demonstrated that, using second-order accurate numerical approximations for the derivatives, the compensation remains fourth-order accurate. Future investigations may be able to specify the effects of the numerical approximation on the simulation results more accurately.

We also gave an example of applying the fourth-order compensation to a finite-element model with larger, sparse matrices. The benefits of the increased accuracy proved to outweigh the computational disadvantages of the non-sparse compensated matrices. This tradeoff could be studied further for even larger systems. We also gave some suggestions on further improving the performance of the implementation of our compensated method, which may also play an advantageous role in this respect.

% Effects of the compensations on efficient performance of matrix operations could also be studied.

% In this paper, we have not discussed the effects of the compensation terms on the shape of the original system matrices. Many numerical solvers exploit these matrices' special structure (e.g., sparseness) for better performance, and breaking this structure through compensation could be detrimental. A rigorous study of the structure of the compensation matrices is thus needed.

So far, we have focused on the Newmark method applied in its classical form, namely, for linear systems under transient loads. Nevertheless, the Newmark method is also widely used on structures with nonlinear internal forces: deriving the distorted equations and investigating the possibility of introducing compensation terms in such cases could be explored in further works.
Compensations of this kind would bring in new derivative terms with new tensorial properties, and an efficient implementation is expected to be an additional challenge.

Additionally, several widely used generalizations of the Newmark method exist. Deriving the distorted equations for such schemes would yield additional insights into the differences and improvements introduced by the generalizations.

The compensation technique introduced here has the potential to improve both the qualitative and quantitative accuracy of numerical simulations performed using the Newmark method. Their greatest advantage is the fact that the improvements can be achieved for the one-time cost of calculating the compensated terms in advance. This means \emph{that no modification of the algorithm implementing the numerical scheme is necessary} -- such a step often poses a significant obstacle in the wider adoption of novel numerical schemes that do require the implementation of new algorithms. We plan to further investigate the possibility of introducing similar compensations for other widely used numerical schemes, such as the generalized-$\alpha$ method showcased here, or other generalizations of the Newmark method mentioned above.

\section*{Acknowledgements}

We thank Bal\'{a}zs T\'{o}th for his valuable suggestions about the scope of our manuscript. 
We are also thankful to Bence Sipos for his insightful analogy between the compensation technique introduced here and the strategy of input preshaping in robotics.

The work was supported by the grant of National Research, Development and Innovation Office -- NKFIH FK 134277, as well as their University Researcher Scholarship Program (EKÖP) which is also supported by the Ministry of Culture and Innovation. We also acknowledge KIFÜ (Governmental Agency for IT Development, Hungary) for awarding us access to the Komondor HPC facility based in Hungary. We also thank the HPC staff at KIFÜ for their technical support.

\section*{Declarations}

\subsection*{Author contributions}
\textbf{Don\'{a}t M. Tak\'{a}cs:} Conceptualization, Methodology, Formal analysis, Software, Validation, Visualization, Writing -- Original Draft.
\textbf{Tam\'{a}s F\"{u}l\"{o}p:} Methodology, Formal analysis, Supervision, Writing -- Review \& Editing.

\subsection*{Conflict of interest}
The authors have no conflict of interest to declare.

\begin{appendices}

\section{}

\subsection{Wolfram Mathematica code for obtaining distorted equations of a general \texorpdfstring{$n$}{n}-dimensional ODE}\label{app:wolfram}
Here, we illustrate the usage of the algorithm on the implicit Euler system applied to a two-dimensional system.

\begin{lstlisting}[language=Mathematica]
y = {y1, y2};
fcn[y1_, y2_] := {y2, -y1}
Phi[h_, {y1j_, y2j_}, {y1jp1_, y2jp1_}] := {y1j, y2j} + h fcn[y1jp1, y2jp1] (* Implicit Euler scheme *)
nn = 5; (* Number of expansion terms *)
fcoe = ConstantArray[0, nn + 1];
diffy = ConstantArray[0, nn + 2];
fcoe[[1]] = fcn @@ y;
For[n = 2, n <= nn, n++,
  modeq = Total@Table[h^j fcoe[[j + 1]], {j, 0, n - 2}];
  diffy[[1]] = y;
  For[i = 1, i <= n, i++,
   diffy[[i + 1]] = D[diffy[[i]], {y}].modeq;
   ];
  ytilde = Total@Table[h^kk diffy[[kk + 1]]/kk!, {kk, 0, n}];
  res = ytilde - (Phi[h, y, ytilde]);
  tay = Sum[h^j/j! (D[res, {h, j}] /. h -> 0), {j, 0, n}];
  fcoe[[n]] = -Coefficient[tay, h, n];
  ];
DVF = Collect[Simplify@Total@Table[h^j fcoe[[j + 1]], {j, 0, nn - 1}], h]
\end{lstlisting}

\subsection{Wolfram Mathematica code for obtaining distorted equations of the Newmark method}\label{app:wolfram2}
\begin{lstlisting}[language=Mathematica]
$Assumptions = Element[Minv, Matrices[{d, d}, Reals, Symmetric[{1, 2}]]] && Element[KK, Matrices[{d, d}, Reals, Symmetric[{1, 2}]]] && Element[CC, Matrices[{d, d}, Reals, Symmetric[{1, 2}]]] && Element[q, Vectors[d, Reals]] && Element[v, Vectors[d, Reals]] && Element[F[tau], Vectors[d, Reals]] && Element[F'[tau], Vectors[d, Reals]] && Element[F''[tau], Vectors[d, Reals]] && Element[h, Reals] && Element[tau, Reals] && Element[beta, Reals] && Element[gamma, Reals];
y = {tau, q, v};
fcn[tau_, q_, v_] := {1, v, -Minv.(CC.v + KK.q - F[tau])};
Phi[h_, {tauj_, qj_, vj_}, {taujp1_, qjp1_, vjp1_}] := With[
   {
    aj = -Minv.(CC.vj + KK.qj - F[tauj]),
    ajp1 = -Minv.(CC.vjp1 + KK.qjp1 - F[taujp1])
    },
   {
    (*taujp1=*)tauj + h,
    (*qjp1=*)qj + h vj + h^2/2 ((1 - 2 beta) aj + 2 beta ajp1),
    (*vjp1=*)vj + h ((1 - gamma) aj + gamma ajp1)
    }
   ];
applyTensorRules[{taudot_, qdot_, vdot_}] := 
 With[{tensorRules = {Dot[1, a___] :> Dot[a], Dot[0, a__] :> 0, 
     Dot[a___, 1] :> Dot[a], Dot[a__, 0] :> 0}}, {taudot, qdot, 
    vdot} //. tensorRules]
nn = 3;(*Number of expansion terms*)
diffy = ConstantArray[0, nn + 2];
fcoe = ConstantArray[0, nn + 1];
fcoe[[1]] = fcn @@ y;
For[n = 2, n <= nn, n++,
  modeq = Total@Table[h^j fcoe[[j + 1]], {j, 0, n - 2}];
  diffy[[1]] = y;
  For[i = 1, i <= n, i++,
   diffy[[i + 1]] = 
     applyTensorRules[
       Inner[Dot, D[TensorExpand@diffy[[i]], {y}], modeq]] /. 
      MatrixPower[a_, k_] :> Dot @@ ConstantArray[a, k];
   ]; 
  ytilde = Total@Table[h^kk diffy[[kk + 1]]/kk!, {kk, 0, n}];
  res = ytilde - (Phi[h, y, ytilde]);
  tay = applyTensorRules[
    Sum[h^j/j! (D[res, {h, j}] /. h -> 0), {j, 0, n}]];
  fcoe[[n]] = -applyTensorRules@Coefficient[tay, h, n];
  ]
MVF = FullSimplify@TensorExpand[Collect[Total@Table[ExpandAll[applyTensorRules[h^j fcoe[[j + 1]]]], {j, 0, nn - 1}], y]] /. MatrixPower[a_, k_] :> Dot @@ ConstantArray[a, k]
\end{lstlisting}

\end{appendices}

\printbibliography

\end{document}